\input ./style/arxiv-general.cfg
\documentclass[aap,MSNbibl,seceqn,dvips]{arximspdf}
\makeatletter
   \@ifpackageloaded{graphicx}{}{\usepackage{graphicx}}
\makeatother

\doi{10.1214/15-AAP1132}
\volume{26}
\issue{3}
\pubyear{2016}
\firstpage{1743}
\lastpage{1773}
\docsubty{FLA}

\makeatletter
\newcommand{\eqref}[1]{(\ref{#1})}
\newcommand{\R}{\mathbb{R}}
\newcommand{\E}{\mathbb{E}}
\newcommand{\PP}{\mathbb{P}}

\newcommand{\D}{\mathbb{D}}

\newcommand{\LLL}{\mathbf{L}}
\newcommand{\al}{\alpha}
\newcommand{\la}{\lambda}
\newcommand{\La}{\Lambda}
\newcommand{\ga}{\gamma}

\newcommand{\si}{\sigma}

\newcommand{\be}{\beta}
\newcommand{\ep}{\varepsilon}

\newcommand{\de}{\delta}
\newcommand{\De}{\Delta}
\newcommand{\om}{\omega}
\newcommand{\Om}{\Omega}
\newcommand{\ze}{\zeta}
\newcommand{\aaa}{\mathcal{A}}
\newcommand{\ba}{\mathcal{B}}
\newcommand{\ea}{\mathcal{E}}
\newcommand{\f}{\mathcal{F}}
\newcommand{\g}{\mathcal{G}}
\newcommand{\h}{\mathcal{H}}
\newcommand{\laa}{\mathcal{L}}
\newcommand{\n}{\mathcal{N}}
\newcommand{\p}{\mathcal{P}}
\newcommand{\s}{\mathcal{S}}
\newcommand{\ua}{\mathcal{U}}
\newcommand{\w}{\mathcal{W}}
\newcommand{\wy}{\widehat{y}}
\newcommand{\BC}{\overline{C}}
\newcommand{\BH}{\overline{H}}
\newcommand{\BY}{\overline{Y}}
\newcommand{\BZ}{\overline{Z}}
\newcommand{\Bf}{\overline{f}}
\newcommand{\Bxi}{\overline{\xi}}
\newcommand{\toop}{\stackrel{\PP}{\longrightarrow}}
\newcommand{\sign}{\operatorname{sign}}
\newtheorem{prop}{Proposition}

\newproclaim{defi}[prop]{Definition}
\newtheorem{lem}[prop]{Lemma}
\newtheorem{theo}[prop]{Theorem}
\newproclaim{rema}[prop]{Remark}
\newproclaim{ex}[prop]{Example}
\newproclaim{ass}{Assumption}
\makeatother

\begin{document}
\begin{frontmatter}

\title{Backward stochastic differential equation driven by a marked point process: An elementary approach with an application to optimal control}
\runtitle{BSDEs driven by a marked point process}

\begin{aug}
\author[A]{\fnms{Fulvia}~\snm{Confortola}\thanksref{M1,t1}\corref{}\ead[label=e1]{fulvia.confortola@polimi.it}},
\author[A]{\fnms{Marco}~\snm{Fuhrman}\thanksref{M1,t1}\ead[label=e2]{marco.fuhrman@polimi.it}}
\and
\author[B]{\fnms{Jean}~\snm{Jacod}\thanksref{M2}\ead[label=e3]{jean.jacod@upmc.fr}}
\runauthor{F. Confortola, M. Fuhrman and J. Jacod}
\affiliation{Politecnico di Milano\thanksmark{M1}    and Universit\'e Pierre et Marie Curie\thanksmark{M2}}
\address[A]{F. Confortola\\
M. Fuhrman\\
Politecnico di Milano\\
Dipartimento di Matematica\\
via Bonardi 9\\
20133 Milano\\
Italy\\
\printead{e1}\\
\phantom{E-mail:\ }\printead*{e2}}
\address[B]{J. Jeacod\\
Institut de Math\'ematiques de Jussieu\\
CNRS-UMR 7586\\
Universit\'e Pierre et Marie Curie\\
4 Place Jussieu\\
75252 Paris-Cedex 5\\
France\\
\printead{e3}}
\end{aug}
\thankstext{t1}{Supported in part  by MIUR-PRIN 2010--2011
 ``Evolution differential problems: deterministic
 and stochastic approaches and their interactions''.}

%
\received{\smonth{7} \syear{2014}}
\revised{\smonth{7} \syear{2015}}

\begin{abstract}
We address a class of
backward stochastic differential equations on a bounded
interval, where the driving noise
is a marked, or multivariate, point process. Assuming that
the jump times are totally inaccessible and a technical
condition holds (see Assumption \ref{aa} below), we
prove existence and uniqueness results under Lipschitz
conditions on the coefficients. Some counter-examples
show that our assumptions are indeed needed.
We use a novel approach that allows reduction to
a (finite or infinite) system of deterministic differential equations,
thus avoiding the use of martingale representation theorems
and allowing potential use of standard numerical methods.
Finally, we apply the main results to solve an optimal control
problem for a marked point process, formulated in
a classical way.
\end{abstract}

\end{frontmatter}

\section{Introduction}\label{s-Intro}

Since the paper  by Pardoux and Peng \cite{PP90}, the topic of backward
stochastic differential equations (BSDE in short) has been in constant
development, due to its utility in finance (see, e.g., El Karoui, Peng and
Quenez \cite{EPQ}), in control theory, and in the theory of nonlinear
PDEs.

The first papers, and most of the subsequent ones, assume that the driving
term is a Brownian motion, but the case of a discontinuous driving process
has also been considered rather early; see, for example, Buckdahn and Pardoux
\cite{BP94}, Tang and Li \cite{tali} and more recently Barles, Buckdahn and
Pardoux \cite{BBP}, Xia \cite{xia}, Becherer \cite{Bech}, Cr\'epey and
Matoussi \cite{CM-2008}, or Carbone, Ferrario and Santacroce
\cite{CaFeSa} among many others.

The case of a driving term which is purely discontinuous has
attracted less attention; see, however, Shen and Elliott \cite{SE-2011}
for the particularly simple ``one-jump'' case, or Cohen and Elliott
\cite{CoEl,CoEl-a} and Cohen and Szpruch \cite{CoSz} for BSDEs
associated to Markov chains. The pure jump case has
certainly less potential applications than the continuous or
continuous-plus-jumps case, but on the other hand it exhibits a~much
simpler structure, which provides original insight on BSDEs.

To illustrate the latter point, in this paper we consider BSDEs driven by
a~marked (or, multivariate) point process. The time horizon is a finite
(nonrandom) time $T$. The point process is nonexplosive, that is, there
are almost surely finitely many points within the interval
$[0,T]$, and it is also \textit{quasi-left continuous}, that is, the jump times
are totally inaccessible: the main examples of this situation are the
Poisson process and the compound Poisson process. We also make the (rather
strong) assumption that the generator is uniformly Lipschitz.

In contrast with most of the literature, in which the martingale
representation theorem and the application of a suitable fixed-point theorem
play a central role, in the setting of point processes it is possible
to solve the equation recursively, by replacing the BSDE by an ordinary
differential equation in between jumps, and match the pre- and post-jump
values at each jump time (such a method has already been used
for a BSDE driven by a Brownian motion plus a Poisson process; see, e.g.,
Kharroubi and Lim \cite{KL-2012}, but then between any two consecutive jumps
one has to solve a genuine BSDE).

Reducing the BSDE to a sequence of ODEs allows us for a very simple solution,
although we still need some elementary  a priori estimates, though,
for establishing the existence when the number of jumps is unbounded.
Apart from the intrinsic interest of a simple method, this might also give
rise to simple numerical ways for solving the equation. Another noticeable
point is that it provides an $\LLL^1$ theory, which is more appropriate for
point processes than the usual $\LLL^2$ theory.

There are two main results about the BSDE: one is when the number of jumps is
bounded, and  then we obtain uniqueness within the class of all possible
solutions. The other is, in the general case, an existence and uniqueness result
within a suitable weighted $\LLL^1$ space. We also state a third important
result, showing how an optimal control problem on a marked process reduces to
solving a BSDE. Existence and uniqueness results for the BSDE are stated
in the case of a scalar equation, but the extension to the vector-valued
case is immediate.

The paper is organized as follows: in Section~\ref{s-1}, we present the setting and the
two main results (as will be seen, the setting is somewhat complicated to
explain, because in the multivariate case there are several distinct but
natural versions for the BSDE). Section~\ref{s-2} is devoted to a few simple
 a priori estimates. In Section~\ref{s-3}, we explain how the BSDE can be reduced
to a sequence of (nonrandom) ODEs, and also exhibit a few counter-examples
when the basic assumptions on the point process are violated. The proof of
the main results is in Section~\ref{s-4}, and in Section~\ref{s-5} the control problem is
considered.

\section{Main results}\label{s-1}
\subsection{The setting}

We have a probability space $(\Om,\f,\PP)$ and a fixed time horizon
$T\in(0,\infty)$, so all processes defined on this space are indexed by $[0,T]$,
and all random times take their values in $[0,T]\cup\{\infty\}$.

This space is endowed with a nonexplosive multivariate point process
(also called marked point process) on $[0,T]\times E$, where $(E,\ea)$ is a
Lusin space: this is a sequence $(S_n,X_n)$ of points, with
distinct times of occurrence $S_n$ and with marks $X_n$, so it can be
viewed as a random measure of the form
\begin{eqnarray}
\label{S1}
&&\mu(dt,dx)=\sum_{n\geq1: S_n\leq T}
\ep_{(S_n,X_n)}(dt,dx),
\end{eqnarray}
where $\ep_{(t,x)}$ denotes the Dirac measure.
Here, the $S_n$'s are $(0,T]\cup\{\infty\}$-valued and the $X_n$'s are
$E$-valued, and $S_1>0$, and $S_n<S_{n+1}$ if $S_n\leq T$, and $S_n\leq
S_{n+1}$ everywhere, and $\Om=\bigcup\{S_n>T\}$. Note that the ``mark'' $X_n$ is
relevant on the
set $\{S_n\leq T\}$ only,  but it is convenient to have it defined on the
whole set $\Om$, and without restriction we may assume that $X_n=\De$ when
$S_n=\infty$, where $\De$ is a~distinguished point in $E$.

We denote by $(\f_t)_{t\geq0}$ the filtration generated by the point process,
which is the smallest filtration for which each $S_n$ is a stopping time and
$X_n$ is $\f_{S_n}$-measurable. As we will see, the special structure of this
filtration plays a fundamental role in all what follows. We let
$\p$ be the predictable $\si$-field on $\Om\times[0,T]$, and for any auxiliary
measurable space $(G,\g)$ a function on the product $\Om\times[0,T]\times G$
which is measurable with respect to $\p\otimes\g$ is called
\textit{predictable}.

We denote by $\nu$ the predictable compensator of the measure $\mu$,
relative to the filtration $(\f_t)$. The
measure $\nu$ admits the disintegration
\begin{eqnarray}
\label{S2}
&& \nu(\om,dt,dx)= dA_t(\om) \phi_{\om,t}(dx),
\end{eqnarray}
where $\phi$ is a transition probability from $(\Om\times[0,T],\p)$ into
$(E,\ea)$, and $A$ is an increasing c\`adl\`ag predictable process
starting at $A_0=0$, which is also the predictable compensator of the
univariate point process
\begin{eqnarray}
\label{S3}
&& N_t = \mu\bigl([0,t]\times E\bigr)~=~\sum
_{n\geq1}1_{\{S_n\leq t\}}.
\end{eqnarray}
Of course, the multivariate point process $\mu$ reduces to the univariate $N$
when $E$ is a singleton.

Unless otherwise specified, the following assumption,
where we set $S_0=0$, will hold throughout.
\renewcommand{\theass}{(A)}
\begin{ass}\label{aa}
(A$_1$)  The process $A$ is continuous (equivalently: the
jump times $S_n$ are totally inaccessible).

(A$_2$)  $\PP(S_{n+1}>T\vert \f_{S_n})>0$ for all $n\ge 0$.
\end{ass}

The first condition amounts to the quasi-left continuity of $N$. We will
briefly examine what happens when (A$_1$) and (A$_2$) fail in Section~\ref{s-3}.

\subsection{The BSDE in the univariate case}

Now, we turn to the BSDE. In addition to the driving point process, the
ingredients are:
\begin{itemize}
\item a \textit{terminal condition} $\xi$, which is always an $\f_T$-measurable
random variable;
\item a \textit{generator} $f$, which is real-valued function depending on $\om$,
on time, possibly on the mark $x$ of the point process, and also in a suitable
way on the solution of the BSDE. In all cases below, the dependence of the
generator upon the solution will be assumed Lipschitz, typically involving
two nonnegative constants $L,L'$, as specified below.
\end{itemize}

We begin with the univariate case, which is
simpler to formulate. In this case, the BSDE takes the form
\begin{eqnarray}
\label{S4}
&& Y_t+\int_{(t,T]}Z_s
\,dN_s=\xi+\int_{(t,T]}f(\cdot,s,Y_{s-},Z_s)
\,dA_s,
\end{eqnarray}
where $f$ is a predictable function on $\Om
\times[0,T]\times\R\times\R$, satisfying
\begin{eqnarray}
\bigl|f\bigl(\om,t,y',z'\bigr)-f(\om,t,y,z)\bigr| &\leq &
L'\bigl|y'-y\bigr|+L\bigl|z'-z\bigr|,
\nonumber
\\[-8pt]
\label{S5}
\\[-8pt]
\nonumber
\int_0^T \bigl|f(t,0,0)\bigr| \,dA_t &<& \infty\qquad
\mbox{a.s.}
\end{eqnarray}

A \textit{solution} is a pair $(Y,Z)$ consisting in an adapted c\`adl\`ag
process $Y$ and a predictable process $Z$ satisfying $\int_0^T|Z_t| \,dA_t
<\infty$ a.s., such that \eqref{S4} holds for all $t\in[0,T]$, outside a
$\PP$-null set [this implicitly supposes that $\int_0^T|f(\cdot,s,Y_s,Z_s)|
 \,dA_s<\infty$ a.s.].

\begin{rema}\label{RS1}  Quite often the BSDE is written, in a slightly
different form, as
\begin{eqnarray}
\label{S6}
&& Y_t+\int_{(t,T]}Z_s
(dN_s-dA_s) =\xi+\int_{(t,T]}f(\cdot,s,x,Y_{s-},Z_s)
\,dA_s.
\end{eqnarray}
Upon a trivial modification of $f$, this is clearly the same as \eqref{S4},
and it explains the integrability restriction on $Z$. The reason underlying
the formulation \eqref{S6} is that it singles out the ``martingale
increment'' $\int_{(t,T]}Z_s (dN_s-dA_s)$.
\end{rema}

\subsection{The BSDE in the multivariate case}

In the multivariate case, the predictable process $Z$ of
\eqref{S4} should be replaced by a predictable function $Z(\om,t,x)$ on
$\Om\times[0,T]\times E$, and this function may enter the generator in
different guises. We start with the most general formulation, and will
single out two special, easier to formulate, cases afterward.

We need some additional notation: we let $\ba(E)$ be the set of all Borel
functions on $E$; if $Z$ is a measurable function on $\Om\times[0,T]\times E$,
we write $Z_{\om,t}(x)=Z(\om,t,x)$, so each $Z_{\om,t}$, often abbreviated
as $Z_t$ or $Z_t(\cdot)$, is an element of $\ba(E)$.

With this notation, the BSDE takes the form
\begin{eqnarray}
&& Y_t+\int_{(t,T]}\int
_EZ(s,x) \mu(ds,dx)
\nonumber
\\[-8pt]
\label{S11}
\\[-8pt]
\nonumber
&&\qquad =\xi+\int_{(t,T]}\int
_Ef\bigl(\cdot,s,x,Y_{s-},Z_s(\cdot)
\bigr) \nu(ds,dx),
\end{eqnarray}
where $f$ is a real-valued function on $\Om
\times[0,T]\times E\times\R \times\ba(E)$, such that
$f\bigl(\om,t,x,y,Z_{\om,t}(\cdot)\bigr)$ is predictable for any
predictable function $Z$ on $\Om\times[0,  T]\times E$, and
\begin{eqnarray}
&& \bigl|f\bigl(\om,t,x,y',\ze\bigr)-f(\om,t,x,y,\ze)\bigr| \leq
L'\bigl|y'-y\bigr|,\nonumber
\\
&& \int_E\bigl|f(\om,t,x,y,\ze)-f\bigl(\om,t,x,y,
\ze'\bigr)\bigr|\phi_{\om,t}(dx)
\nonumber
\\[-8pt]
\label{S12}
\\[-8pt]
\nonumber
&&\qquad\leq   L\int_E\bigl|
\ze'(x)-\ze(x)\bigr| \phi_{\om,t}(dx),
\\
&& \int_0^T\!\int_E\bigl|f(t,x,0,0)\bigr|
\nu(dt,dx)< \infty\qquad \mbox{a.s.} \nonumber
\end{eqnarray}
[in the expression $f(t,x,0,0)$, the last ``$0$'' stands for the function
in $\ba(E)$ which vanishes identically].

A \textit{solution} is a pair $(Y,Z)$ consisting in an adapted c\`adl\`ag
process $Y$ and a predictable function $Z$ on $\Om\times[0,T]\times E$
satisfying $\int_0^T\!\int_E|Z(t,x)| \nu(ds,dx)<\infty$ a.s., such
that \eqref{S11} holds for all $t\in[0,T]$, outside a $\PP$-null set.

The measurability condition imposed on the generator is somewhat awkward,
and probably difficult to check in general. However, it is satisfied in
the following two types of equations.

\begin{longlist}[\textit{Type II equation}:]
\item[\textit{Type I equation}:]  This is the simplest one to state, and
it takes the form
\begin{eqnarray}
&& Y_t+\int_{(t,T]}\int
_EZ(s,x) \mu(ds,dx)
\nonumber
\\[-8pt]
\label{S7}
\\[-8pt]
\nonumber
&&\qquad=\xi+\int_{(t,T]}\int
_Ef_I\bigl(\cdot,s,x,Y_{s-},Z(s,x)
\bigr) \nu(ds,dx),
\end{eqnarray}
where $f_I$ is a predictable
function on $\Om\times[0,T]\times E\times\R\times\R$, satisfying
\begin{eqnarray}
\bigl|f_I\bigl(\om,t,x,y',z'
\bigr)-f_I(\om,t,x,y,z)\bigr| &\leq &  L'\bigl|y'-y\bigr|+L\bigl|z'-z\bigr|,
\nonumber
\\[-8pt]
\label{S8}
\\[-8pt]
\nonumber
\int_0^T\!\int_E\bigl|f_I(t,x,0,0)\bigr|
\nu(dt,dx)&<&\infty\qquad \mbox{a.s.}
\end{eqnarray}

That \eqref{S7} is a special case of \eqref{S11} is obvious; we simply
have to take for $f$ the function on $\Om\times[0,T]\times E\times\R
\times\ba(E)$ defined by
\begin{eqnarray}
\label{S13}
&& f(\om,s,x,y,\ze)=f_I\bigl(\om,s,x,y,\ze(x)\bigr),
\end{eqnarray}
and \eqref{S8} for $f_I$ yields \eqref{S12} for $f$.

\item[\textit{Type II equations}:]  The BSDE \eqref{S7} cannot in general be used
as a tool for solving control problems driven by a multivariate
point process, whereas this is one of the main motivations for introducing
them. We rather need the following formulation:
\begin{eqnarray}
\label{S9}\quad
&& Y_t+\int_{(t,T]}\int
_EZ(s,x) \mu(ds,dx) =\xi+\int_{(t,T]}f_{\mathrm{II}}(\cdot,s,Y_{s-},
\eta_sZ_s) \,dA_s,
\end{eqnarray}
where, recalling that $\phi_{\om,t}$ are the measures occurring in
\eqref{S2} and $Z_{\om,t}(x)=Z(\om,t,x)$,
\begin{eqnarray}
&& \mbox{$\eta_{\om,t}$ is a real-valued
map on $\ba(E)$},\nonumber
\\
&&\qquad\mbox{with }\bigl|\eta_{\om,t}\ze-\eta_{\om,t}
\ze'\bigr|\leq \int_E\bigl|\ze'(v)-
\ze(v)\bigr| \phi_{\om,t}(dv),\nonumber
\\
\label{S10}&&
\mbox{$Z$ predictable on }\Om\times[0,T]\times E \quad\Rightarrow\\
&&\qquad\mbox{the
process $(\om,t)\mapsto \eta_{\om,t}Z_{\om,t}$ is predictable,}
\nonumber\\
&&\mbox{$f_{\mathrm{II}}$ is a function satisfying
\eqref{S5}.}\nonumber
\end{eqnarray}

Again, \eqref{S9} reduces to \eqref{S11} upon taking
\begin{eqnarray}
\label{S14}
&& f(\om,s,x,y,\ze)=f_{\mathrm{II}}(\om,s,y,\eta_{\om,s}
\ze),
\end{eqnarray}
and \eqref{S5} for $f_{\mathrm{II}}$ plus \eqref{S10} for $\eta_{\om,t}$
yield \eqref{S12} for $f$. As we will see in Section~\ref{s-5}, this type
of equation is well suited to control problem.
\end{longlist}

In the univariate case, all three formulations \eqref{S11}, \eqref{S7}
and \eqref{S9} coincide with~\eqref{S4}.

Finally, we describe another notion of a solution, starting with the
following remark: we can of course rewrite \eqref{S11} as follows:
\begin{eqnarray}
&& Y_t+\sum_{n\geq1}Z(S_n,X_n)
1_{\{t<S_n\leq T\}}
\nonumber
\\[-8pt]
\label{S15}
\\[-8pt]
\nonumber
&&\qquad=\xi+\int_{(t,T]}\int_Ef
\bigl(s,x,Y_{s-},Z_s(\cdot)\bigr) \nu(ds,dx).
\end{eqnarray}
Since $A$ is continuous, \eqref{S15} yields, outside a $\PP$-null set,
\begin{eqnarray}
&& \De Y_{S_n}=Z(S_n,X_n)\qquad\mbox{if
$S_n\leq T$ and $n\geq1$},
\nonumber
\\[-8pt]
\label{S16}
\\[-8pt]
\nonumber
&& Y\mbox{ is continuous outside $
\{S_1,\ldots,S_n,\ldots\}$}.
\end{eqnarray}
In other words, $Y$ completely determines the predictable function
$Z$ outside a~null set with respect to the measure
 $\PP(d\om)\mu(\om,dt,dx)$, hence also
outside a~$\PP(d\om)\nu(\om,dt,dx)$-null set. Equivalently, if $(Y,Z)$ is a
solution and $Z'$ is another predictable function, then $(Y,Z')$ being another
solution is the same as having $Z'=Z$ outside a $\PP(d\om)\mu(\om,dt,
dx)$-null set, and the same as having $Z'=Z$ outside a $\PP(d\om)
\nu(\om,dt,dx)$-null set.

Therefore, another way of looking at equation \eqref{S11} is as follows:
a \textit{solution} is an adapted c\`adl\`ag process $Y$ for which there exists
a predictable function $Z$ satisfying
\[
\int_0^T\!\int_E\bigl|Z(s,x)\bigr| \nu(ds,dx)
<\infty \qquad\mbox{a.s.},
\]
such that the pair $(Y,Z)$ satisfies \eqref{S11} for all
$t\in[0,T]$, outside a $\PP$-null set. Then \textit{uniqueness} of the solution
means that, for any two solutions $Y$ and $Y'$ we have $Y_t=Y'_t$ for all
$t\in[0,T]$, outside a $\PP$-null set.

\subsection{Statement of the main results}
We have two main results. The first one is when the point
process has at most $M$ points, for a nonrandom integer $M$, that is,
\begin{eqnarray}
\label{S20}
&& \PP(S_{M+1}=\infty)=1.
\end{eqnarray}

\begin{theo}\label{TS1}
Assume \ref{aa} and \eqref{S20}. The solution $Y$ of
\eqref{S11}, if it exists, is unique up to
null sets. Moreover, if the variable $A_T$ is bounded, and if
\begin{eqnarray}
\label{S21}
&& \E\bigl(|\xi|\bigr)<\infty,\qquad \E \biggl(\int_0^T\!
\int_E\bigl|f(s,x,0,0)\bigr| \nu(ds,dx) \biggr)<\infty,
\end{eqnarray}
the solution exists and satisfies $\E (\int_0^T|Y_t| \,dA_t )<\infty$
and $\E (\int_0^T\int_E|Z(t,x)| \nu(dt,\break dx) )<\infty$.
\end{theo}

The existence result above is ``almost'' a special case of the next
theorem. In contrast, the uniqueness within the class of \textit{all} possible
solutions is specific to the situation \eqref{S20}. When this fails,
uniqueness holds only within smaller subclasses, which we now describe.
For any $\al>0$ and $\be\geq0$, we set
\begin{eqnarray}
\laa_{\al,\be}^1 &=& \mbox{the set
of all pairs $(Y,Z)$ with $Y$ c\`adl\`ag adapted and}\nonumber\\
\label{S22}
&&{}\mbox{$Z$ predictable, satisfying,}
\\
\nonumber
\bigl\|(Y,Z)\bigr\|_{\al,\be} &:=&  \E \biggl(\int_0^T\!
\int_E \bigl(|Y_t|+\bigl|Z(t,x)\bigr| \bigr)
e^{\be
A_t} \al^{N_t} \nu(dt,dx) \biggr)<\infty.
\end{eqnarray}
The space $\laa^1_{\al,\be}$ decreases when $\al$ and/or $\be$ increases.

\begin{theo}\label{TS2}
Assume~\ref{aa}.
\begin{longlist}[(a)]
\item[(a)]  If
\begin{eqnarray}
\E \bigl(e^{\be A_T} \al^{N_T} |\xi| \bigr) &<& \infty,
\nonumber
\\[-8pt]
\label{S23}
\\[-8pt]
\nonumber
\E \biggl(\int_0^T\!\int_E
\al^{N_s} e^{\be A_s} \bigl|f(s,x,0,0)\bigr| \nu(ds,dx) \biggr)&<& \infty,
\end{eqnarray}
for some $\al>L$ and $\be>1+\al+L'$, where $L,L'$ are the
constants occurring in~\eqref{S12},
then \eqref{S11} admits one and only one (up to null sets) solution $(Y,Z)$
belonging to $\laa^1_{\al,\be}$.

\item[(b)] When moreover the variable $A_T$ is bounded, the conditions
\begin{eqnarray}
\label{S24}
&& \quad\E \bigl(|\xi|^{1+\ep} \bigr)<\infty,\qquad \E \biggl( \biggl(
\int_0^T\!\int_E\bigl|f(s,x,0,0)\bigr|
\nu(ds,dx) \biggr)^{1+\ep} \biggr)<\infty
\end{eqnarray}
for some $\ep>0$ imply \eqref{S23} for all $\be\geq0$ and $\al>0$, hence
\eqref{S11} admits one and only one (up to null sets) solution $(Y,Z)$
belonging to $\bigcup_{\al>L, \be>1+\al+L'} \laa^1_{\al,\be}$, and this
solution also belongs to $\bigcap_{\al>0,\be\geq0} \laa^1_{\al,\be}$.
\end{longlist}
\end{theo}

The claim (b) is interesting, because it covers the most usual situation where
$\mu$ is a Poisson random measure (so that $A_t=\la t$ for some constant
$\la>0$). Note that, even in this case, we do not know whether \eqref{S11}
admits other solutions, which are not in\vspace*{1pt} $\bigcup_{\al>L, \be>1+\al+L'}
\laa^1_{\al,\be}$.

We note that if we apply Theorem~\ref{TS2} with the assumptions of
Theorem~\ref{TS1}, namely $A_T \le K$ and $N_T \le M$, condition \eqref{S23}
is equivalent to \eqref{S21} since the exponential factors
are bounded. In this sense, Theorem~\ref{TS1} is a special case of Theorem~\ref{TS2}, except that in the latter theorem uniqueness is guaranteed only
within the smaller class $\laa_{\al,\be}^1$. The occurrence of
exponential weights in the definition of the norm in this space is due
to the fact that we are dealing with BSDEs driven by a general random
compensator $\nu(\om,dt,dx)=dA_t(\om) \phi_{\om,t}(dx)$, where\vspace*{1pt} $A$ is an
increasing  but not necessarily bounded predictable processes.
The same happens in the
$\LLL^2$ theory for BSDEs associated to  marked point processes (see
\cite{CoFu-m,xia}) and for BSDEs driven by a general c\`{a}dl\`{a}g
martingale (see \cite{EH}). On the other hand, in case of compensators
absolutely continuous with respect to a deterministic measure, \cite{tali,Bech,CM-2008},  a standard  $\LLL^2$  theory holds (the norm
reduces to a simpler form, not involving exponentials of stochastic
processes).

\section{A priori estimates}\label{s-2}

In this section, we provide some  a priori estimates for the solutions
of equation \eqref{S11}. Without special mention, Assumption (A$_1$) is
assumed throughout.

\begin{lem}\label{LP1} Let $\al>0$ and $\be\in\R$. If $(Y,Z)$
is a solution of \eqref{S11} we have almost surely
\begin{eqnarray}
&& |Y_t|e^{\be A_t}
\al^{N_t}+\int_t^T\!\int
_E \bigl(\al\bigl|Y_{s-}+Z(s,x)\bigr| -|Y_{s-}|
\bigr) e^{\be A_s} \al^{N_{s-}} \mu(ds,dx)
\nonumber\\
\label{P1}
&&\quad{}+\be\int
_t^T|Y_s| e^{\be A_s}
\al^{N_s} \,dA_s
\\
&&\qquad =|\xi|^p e^{\be A_T} \al^{N_T}+\int
_t^T\!\int_E
\sign(Y_s) f\bigl(s,x,Y_s,Z_s(\cdot)\bigr)
e^{\be A_s} \al^{N_s} \nu(ds,dx).\nonumber
\end{eqnarray}
\end{lem}

\begin{pf}
Letting $U_t$ and $V_t$ be the left-hand and right-hand sides of
\eqref{P1}, and since these processes are c\`adl\`ag, and continuous
outside the $S_n$'s, and $U_T=V_T$, it suffices to check that outside a null
set we have $\De U_{S_n}=
\De V_{S_n}$ and also $U_t-U_s=V_t-V_s$ if $S_n\leq t<s<S_{n+1}\wedge T$, for
all $n\geq0$. The first property is obvious because $\De Y_{S_n}=Z(S_n,X_n)$
a.s. and $A$ is continuous. The second property follows from
$Y_t-Y_s=\int_t^s\!\int_Ef(v,x,Y_v,Z_v(\cdot)) \nu(dv,dx)$,
implying $|Y_t|-|Y_s|=\int_t^s\!\int_E\sign(Y_v) f(v,x,Y_v,
Z_v(\cdot)) \nu(dv,dx)$
and $\al^{N_v}=\al^{N_t}$ for all $v\in[t,s]$, plus a standard change
of variables formula.
\end{pf}

For any $\al>0$ and $\be\geq0$, and with any measurable process $Y$ and
measurable function $Z$ on $\Om\times[0,T]\times E$ we set for
$0\leq t<s\leq T$
\begin{eqnarray}
\label{P2}
&& \w^{\al,\be}_{(t,s]}(Y,Z)=\int_t^s
\!\int_E \bigl(|Y_v|+\bigl|Z(v,x)\bigr| \bigr)
e^{\be
A_v} \al^{N_{v}} \nu(dv,dx),
\end{eqnarray}
so with the notation \eqref{S22} we have $\|(Y,Z)\|_{\al,\be}=
\E(\w^{\al,\be}_{(0,T]}(Y,Z))$. Below, $L$ and $L'$ are as in \eqref{S12}.

\begin{lem}\label{LP2} Let $\al>L$ and $\be>1+\al+L'$. There is a constant
$C$ only depending on $(\al,\be,L,L')$, such
that\vspace*{1pt} any pair $(Y,Z)$ in $\laa^1_{\al,\be}$ which solves \eqref{S11}
satisfies, for any stopping time $S$ with $S\leq T$ and outside a null set,
\begin{eqnarray}
&& \quad|Y_ S|e^{\be A_S} \al^{N_S}
\nonumber
\\[-8pt]
\label{P3}
\\[-8pt]
\nonumber
&&\quad\qquad\leq  \E
\biggl(|\xi|e^{\be A_T} \al^{N_T}+ \int_S^T
\!\int_E\bigl|f(s,x,0,0)\bigr| e^{\be A_s} \al^{N_s}
\nu(ds,dx)\Big\vert \f_S \biggr),
\\
&&\quad\E\bigl(\w^{\al,\be}_{(S,T]}(Y,Z)\vert
\f_S\bigr)
\nonumber
\\[-8pt]
\label{P4}
\\[-8pt]
\nonumber
&&\quad\qquad\leq   C \E \biggl(|\xi|e^{\be A_T} \al^{N_T}+
\int_S^T\!\int_E
\bigl|f(s,x,0,0)\bigr| e^{\be A_s} \al^{N_s} \nu(ds,dx)\Big\vert
\f_S \biggr).
\end{eqnarray}
\end{lem}

\begin{pf}
We have
$\al|Y_{s-}+Z(s,x)|-|Y_{s-}|\geq\al|Z(s,x)|-(1+\al)|Y_{s-}|$, hence~\eqref{P1}, and the Lipschitz condition \eqref{S12} plus the fact that
$\phi_{t,\om}(E)=1$ yield
almost surely
\begin{eqnarray}
&& |Y_S|e^{\be A_S}
\al^{N_S}+\al\int_S^T\!\int
_E\bigl|Z(s,x)\bigr| e^{\be A_s} \al^{N_{s-}} \mu(ds,dx) +
\be\int_S^T|Y_s|e^{\be A_s}
\al^{N_s} \,dA_s\nonumber
\\
\label{P5} &&\qquad\leq|\xi|e^{\be A_T} \al^{N_T}+ (1+\al)\int_S^T|Y_{s-}|
e^{\be A_s} \al^{N_{s-}} \,dN_s
\\
&&\qquad\quad
{}+\int_S^T\!\int_E
\bigl(\bigl|f(s,x,0,0)\bigr|+L'|Y_s|+L\bigl|Z(s,x)\bigr| \bigr)
e^{\be A_s} \al^{N_s} \nu(ds,dx).
\nonumber
\end{eqnarray}

Since $\E(\int_S^T\!\int_E\psi(s,x) \mu(ds,dx)\vert \f_S)=
\E(\int_S^T\!\int_E\psi(s,x) \nu(ds,dx)\vert \f_S)$ for any nonnegative
predictable function $\psi$,
taking the $\f_S$-conditional expectation in \eqref{P5} yields
\begin{eqnarray*}
&& |Y_S| e^{\be A_S} \al^{N_S}+\E \biggl(\int_S^T\!\int_E \bigl(\al\bigl|Z(s,x)\bigr|
+\be|Y_s| \bigr) e^{\be A_s} \al^{N_{s}} \nu(ds,dx)\Big\vert \f_S \biggr)\\
 &&\qquad    \leq \E \bigl(|\xi|e^{\be A_T} \al^{N_T} \bigr)
  \\
  &&\qquad\quad{}+
\E \biggl( \int_S^T  \int_E \bigl(\bigl|f(s,x,0,0)\bigr|+\bigl(1+\al+L'\bigr)|Y_{s-}|+L\bigl|Z(s,x)\bigr| \bigr)
\\
&&\qquad\quad {}\times e^{\be A_s} \al^{N_s} \nu(ds,dx)\Big\vert \f_S \biggr).
\end{eqnarray*}
When $\E(\w^{\al,\be}_{(0,T]}(Y,Z))<\infty$, this implies almost surely
\begin{eqnarray*}
&& |Y_S| e^{\be A_S} \al^{N_S}+\E \biggl(\int_S^T\!\int_E \bigl(\bigl(\be-1-\al-L'\bigr)|Y_s|
+(\al-L)\bigl|Z(s,x)\bigr| \bigr)\\
&&\quad{}\times e^{\be A_s} \al^{N_{s}}\nu(ds,dx)\Big\vert \f_S \biggr)\\
&&\qquad
\leq \E \biggl(|\xi|e^{\be A_T}\al^{N_T}+
\int_S^T\!\int_E\bigl|f(s,x,0,0)\bigr| e^{\be A_s} \al^{N_s} \nu(ds,dx)\Big\vert \f_S \biggr),
\end{eqnarray*}
giving us both \eqref{P3} and \eqref{P4}.
\end{pf}

\begin{lem}\label{LP3} Let $\al>L$ and $\be>1+\al+L'$. If
$(Y,Z)$ is a solution of \eqref{S11} and $(Y',Z')$ is a
solution of the same equation with the same generator $f$ and another
terminal condition $\xi'$, both pairs
$(Y,Z)$ and $(Y',Z')$ being in $\laa^1_{\al,\be}$, we have for any stopping
time $S$ with $S\leq T$ and outside a null set
\begin{eqnarray}
\label{P6}
 \bigl|Y'_S-Y_S\bigr|e^{\be A_S}
\al^{N_S} & \leq &  \E \bigl(\bigl|\xi'-\xi\bigr|e^{\be A_T}
\al^{N_T} \vert \f_S \bigr),
\\
\label{P7}
\E\bigl(\w^{\al,\be}_{(0,T]}\bigl(Y'-Y,Z'-Z
\bigr)\bigr) &\leq &  C \E \bigl(\bigl|\xi'-\xi\bigr|e^{\be A_T}
\al^{N_T} \bigr).
\end{eqnarray}
In particular, \eqref{S11} admits, up to null sets, at most one solution
$(Y,Z)$ belonging to $\laa^1_{\al,\be}$.
\end{lem}

\begin{pf}
Set [with $\ze$ arbitrary in $\ba(E)$, and recalling the
notation $Z_{\om,t}(x)=Z(\om,t,x)$]
\begin{eqnarray*}
&& \BY = Y'-Y,\qquad \BZ=Z'-Z,\qquad \Bxi=\xi'-\xi,\\
&& \Bf(\om,s,x,y,\ze)\\
&&\qquad =
f\bigl(\om,s,x,Y_{s-}(\om)+y,Z_{\om,s}(\cdot)+\ze\bigr)-f\bigl(\om,s,x,Y_{s-}(\om),
Z_{\om,s}(\cdot)\bigr).
\end{eqnarray*}
Then $\Bf$ is satisfies \eqref{S12} with the same constants $L,L'$, and
also $\Bf(s,x,0,0)=0$, and clearly $(\BY,\BZ)$ belongs to $\laa^1_{\al,\be}$
and satisfies \eqref{S11} with the generator $\Bf$ and the terminal condition
$\Bxi$. Hence, \eqref{P6} and \eqref{P7} are exactly \eqref{P3} and
\eqref{P4} written for $(\BY,\BZ)$.

Finally, the last claim follows by taking $\xi'=\xi$.
\end{pf}

\section{The structure of the solutions}\label{s-3}

In this section, we show how it is possible to reduce the problem of solving
equation \eqref{S11} to  solving a sequence of ordinary differential
equations. This reduction needs a number of rather awkward notation, but
it certainly has interest in its own sake. Except in the last subsection,
devoted to some counter-examples, we assume \ref{aa}. We stress that both A$_1$
and A$_2$ are crucial here, in particular to characterize the
$\f_{S_n}$-conditional law of $(S_{n+1},X_{n+1})$ and the compensator $\nu$
of $\mu$.

\subsection{Some basic facts}\label{ss-3-1}

Recall that $(S_n,X_n)$ takes its values in the set $\s=([0,T]\times E)
\cup\{(\infty,\De)\}$.
For any integer $n\geq0$, we let $H_n$ be the subset of $\s^{n+1}$
consisting in all $D=((t_0,x_0),\ldots,(t_n,x_n))$ satisfying
\begin{eqnarray*}
&& t_0 = 0, \qquad x_0=\De, \qquad t_{j+1}\geq t_j,\qquad
t_j\leq T \\
&& \qquad\Rightarrow \quad t_{j+1}>t_j,\qquad  t_j>T \\
&& \qquad\Rightarrow\quad
(t_j,x_j)=(\infty,\De).
\end{eqnarray*}
We set $D^{\max}=t_n$ and endow $H_n$ with its Borel $\si$-field $\h_n$.
We set $S_0=0$ and $X_0=\De$, so
\begin{eqnarray}
\label{B2}
&& D_n=\bigl((S_0,X_0),
\ldots,(S_n,X_n)\bigr)
\end{eqnarray}
is a random element with values in $H_n$, whose law is denoted as $\La_n$
[a probability measure on $(H_n,\h_n)$].

The filtration $(\f_t)$ generated by the point process $\mu$ has a very
special structure, which reflects on adapted or predictable processes, and
below we explain some of these properties; see \cite{J74} for more details.
They might look complicated at
first glance, but they indeed allow us to replace random elements by
deterministic functions of all the $D_n$'s.
\begin{longlist}[(a)]
\item[(a)]  \textit{The variable} $\xi$:  Since $\xi$ is $\f_T$-measurable,
for each $n\geq0$ there is an $\h_n$-measurable map $D\mapsto u^n_D$ on $H_n$
with
\begin{eqnarray}
&& D^{\max}=\infty \quad\Rightarrow \quad u^n_D=0,
\nonumber
\\[-8pt]
\label{B3}
\\[-8pt]
\nonumber
&& S_n(\om)\leq T<S_{n+1}(\om)\quad\Rightarrow\quad\xi(
\om)=u_{D_n(\om)}^n.
\end{eqnarray}

\item[(b)] \textit{Adapted c\`adl\`ag processes}:  A c\`adl\`ag process $Y$, which
further is continuous outside the times $S_n$, is adapted if and only if for
each $n\ge 0$ there is a Borel function $y^n=y^n_D(t)$ on $H_n\times[0,T]$ such
that
\begin{eqnarray}
&& D^{\max}=\infty \quad\Rightarrow \quad y^n_D(t)=0,\nonumber
\\
\label{B4}
&& t\mapsto y^n_D(t)\mbox{ is continuous on $[0,T]$ and
constant on $\bigl[0,T\wedge D^{\max}\bigr]$},
\\
&& S_n(\om)\leq t< S_{n+1}(\om),\qquad t\leq T\quad
\Rightarrow \quad Y_t(\om)=y_{D_n(\om)}^n(t),
\nonumber
\end{eqnarray}
and we express this as $Y\equiv(y^n)$.

\item[(c)]  \textit{Predictable functions}:  A function $Z$ on $\Om\times[0,T]\times E$
is predictable if and only if for each $n\ge 0$ there is a Borel function
$z^n=z^n_D(t,x)$ on $H_n\times[0,T]\times E$ such that
\begin{eqnarray}
 D^{\max} &=& \infty\quad\Rightarrow\quad z^n_D(t,x)=0,
\nonumber
\\[-8pt]
\label{B5}
\\[-8pt]
\nonumber
S_n(\om) &<& t\leq S_{n+1}(\om)\wedge T\quad\Rightarrow\quad Z(
\om,t,x)=z_{D_n(\om)}^n(t,x).
\end{eqnarray}
We express this as $Z\equiv(z^n)$, and also write $z^n_{D,t}$ for the
function $z^n_{D,t}(x)=z^n_D(t,x)$ on $E$.

\item[(d)] \textit{The $\f_{S_n}$-conditional law of $(S_{n+1},X_{n+1})$}:
This conditional law takes the form $G^n_{D_n}$, where $G^n_D(dt,dx)$ is
a transition probability from $H_n$ into $[0,\infty]\times E$, and upon
using \ref{aa} we may further assume the following structure on $G^n_D$, where
$\phi^n_{D,t}(dx)$ is a transition probability from $H_n\times[0,\infty]$
into $E$:
\begin{eqnarray}
&& G^n_D(dt,dx)=G'^n_D(dt)
\phi^n_{D,t}(dx)\qquad \mbox{where }G'^n_D(dt)=G^n_D(dt,E),\nonumber
\\
&& G'^n_D\bigl((T,\infty)\bigr)=0, \qquad t>T
\quad\Rightarrow\quad\phi^n_{D,t}(dx)=\ep_\De(dx),
\nonumber\\
\label{B6}
&& t\mapsto g^n_{D}(t):=G'^n_D\bigl((t,
\infty]\bigr)\mbox{ is continuous }\bigl(\mbox{by (A$_1$)}\bigr),
\\
&& g^n_D(T)>0\qquad
\bigl(\mbox{by }(\mathrm{A}_2)\bigr),
\nonumber\\
&& D^{\max}<\infty\quad\Rightarrow\quad g^n_D
\bigl(D^{\max}\bigr)=1.
\nonumber
\end{eqnarray}
The last property $D^{\max}<\infty~\Rightarrow~g^n_D(D^{\max})=1$,
which plays an important role later, simply expresses the fact that
$S_{n+1} > S_n$ if $S_n < \infty$.

\item[(e)]  \textit{The compensator $\nu$ of $\mu$}:  The following gives us versions
of $\nu$ and $A$ and $\phi_{\om,t}$ in \eqref{S2}:
\begin{eqnarray}
\nu(\om;dt,dx) &=& \sum
_{n=0}^\infty \nu^n_{D_n(\om)}(dt,dx)
1_{\{S_n<t\leq S_{n+1}\wedge T\}},
\nonumber\\
\nu^n_D(dt,dx)& =&\frac{1}{g^n_{D}(t)}
G^n_{D}(dt,dx),
\nonumber
\\[-8pt]
\label{B7}
\\[-8pt]
\nonumber
S_n(\om) &<& t\leq S_{n+1}(\om)\quad\Rightarrow\quad
\phi_{\om,t}=\phi^n_{D_n(\om),t},
\\
A_t(\om)&=& \sum_{n=0}^\infty
a_{D_n(\om)}^n\bigl(t\wedge S_{n+1}(\om)\bigr),\qquad
a^n_D(t)=-\log g^n_D(t),
\nonumber
\end{eqnarray}
hence\vspace*{1pt} $a^n_D(t)=0$ for $t\leq D^{\max}$, and $a^n_D(T)<\infty$.

\item[(f)] \textit{The generator}:  Recall that we are interested in equation
\eqref{S11}, so by \eqref{S12} the generator $f$ has a nice predictability
property only after plugging in a predictable function $Z$. This implies
that, for any $n\geq0$, and if $z^n=z^n_D(t,x)$ is as in (c) above, one has
a Borel function $f\{z^n\}^{n}=f\{z^n\}^{n}_D(t,x,y,w)$ on
$H_n\times[0,T]\times E\times\R\times\R$, such that (with $t\leq T$ below)
\begin{eqnarray}
&& D^{\max} = \infty\quad\Rightarrow\quad f\bigl
\{z^n\bigr\}^n_{D}(t,x,y)=0,
\nonumber
\\
\label{B8}
&& S_n(\om)< t\leq S_{n+1}(\om),\qquad\ze(x)=w+z^n_{D_n(\om)}(t,x)\\
&&\qquad\Rightarrow\quad f(\om,t,x,y,\ze)=f\bigl\{z^n\bigr\}^n_{D_n(\om)}(t,x,y,w).\nonumber
\end{eqnarray}
Moreover, the last two conditions in \eqref{S12} imply that one can take a
version which satisfies identically (where $z^n$ and $z'^n$ are two terms
as in (c), and $f\{0\}^n_D$ below is $f\{z^n\}^n_D$
for $z^n_D(t,x)\equiv0$)
\begin{eqnarray}
&&  \bigl|f\bigl\{z^n\bigr\}^n_D\bigl(t,x,y',w'\bigr)-f
\bigl\{z'^n\bigr\}^n_D(t,x,y,w)\bigr|\nonumber
\\
\label{B9}
&&\qquad\leq L'\bigl|y'-y\bigr|+L\bigl|w'-w\bigr| +L\int
_E\bigl|z'^n_D(t,v)-z^n_D(t,v)\bigr|
\phi^n_{D,t}(dv)
\\
\nonumber
&& \int_0^T\bigl|f\{0\}^n_D(t,x,0,0)\bigr|
\nu^n_D(dt,dx)<\infty.
\end{eqnarray}
\end{longlist}

\subsection{Reduction to ordinary differential equations}\label{ss-3-2}

By virtue of \eqref{S16}, if $Y\equiv(y^n)$ is a
solution of \eqref{S11}, we can, and always will, take for the associated
process $Z\equiv(z^n)$ the one defined for $t\in[0,T]$ by
\begin{eqnarray}
\label{B1}
&& z^n_{D}(t,x)=y_{D\cup\{(t,x)\}}^{n+1}(t)
1_{\{t>D^{\max}\}} -y_{D}^{n}(t),
\end{eqnarray}
because\vspace*{1.5pt} $Y_{S_{n+1}}=y^{n+1}_{D_n\cup\{(S_{n+1},X_{n+1})\}}(S_{n+1})$ and
$Y_{S_{n+1}-}=y^{n}_{D_{n}}(S_{n+1})$, when $S_{n+1}\leq T$. We will in
fact write the above in another form, suitable for plugging into the
generator $f$, as represented by \eqref{B8}. Namely, we set
\begin{eqnarray}
 \wy^{n+1}&=&\bigl(\wy^{n+1}_D(t,x):
(D,t,x)\in H_n\times[0,T]\times E\bigr):
\nonumber
\\[-8pt]
\label{B10}
\\[-8pt]
\nonumber
&&{} \wy^{n+1}_D(t,x)=y_{D\cup\{(t,x)\}}^{n+1}(t)
1_{\{t>D^{\max}\}}.
\end{eqnarray}
Then we take $Z\equiv(z^n)$ as follows:
\begin{eqnarray}
\label{B11}
&& z^n_D(t,x)=\wy^{n+1}_D(t,x)-y^n_D(t),
\end{eqnarray}
and it follows that
\begin{eqnarray}
&& S_n(\om)<t\leq S_{n+1}(\om)
\nonumber
\\
&&\label{B12}
\qquad\Rightarrow\quad f
\bigl(\om,t,x,Y_{t-},Z_t(\cdot)\bigr)
\\
&& \qquad\qquad\qquad= f\bigl\{
\wy^{n+1}\bigr\}^n_{D_n(\om)}\bigl(t,x,y^n_{D_n(\om)}(t),-y^n_{D_n(\om)}(t)
\bigr).\nonumber
\end{eqnarray}

The following lemma is a key point for our analysis.

\begin{lem}\label{LB1} A c\`adl\`ag adapted process $Y\equiv(y^n)$ solves
\eqref{S11} if and only if for $\PP$-almost all $\om$ and all
$n\geq0$ we have
\begin{eqnarray}
&& y^n_{D_n(\om)}(t)\nonumber\\
&&\qquad=u^n_{D_n(\om)}
\nonumber
\\[-8pt]
\label{B13}
\\[-8pt]
\nonumber
&&\qquad\quad{}+
\int_t^T\!\int_Ef\bigl\{
\wy^{n+1}\bigr\}^n_{D_n(\om)}\bigl(s,x,y^n_{D_n(\om)}(s),
-y^n_{D_n(\om)}(s)\bigr) \nu^n_{D_n(\om)}(ds,dx),\\
\eqntext{\displaystyle t\in[0,T].}
\end{eqnarray}
If further \eqref{S20} holds, then $Y\equiv(y^n)$ is a solution if and
only if for $\PP$-almost all $\om$ we have \eqref{B13} for all
$n=0,\ldots,M-1$ and
\begin{eqnarray}
\label{B14}
&& t\in[0,T]\quad\Rightarrow\quad y^M_{D_M(\om)}(t)=u^M_{D_M(\om)}=
\xi(\om).
\end{eqnarray}
\end{lem}

\begin{pf}
Considering the restriction of the BSDE to each interval
$[S_n,\break S_{n+1})\cap[0,T]$  and recalling \eqref{S16}, we see that
$Y$ is a solution if and
only if, outside some null set $\n$, we have for $n\geq0$
\begin{eqnarray*}
S_{n} &\leq &  t<S_{n+1}\leq T \quad\Rightarrow\quad
Y_t=Y_{S_{n+1}-}+\int_t^{S_{n+1}}\!\int_Ef\bigl(s,x,Y_s,Z_s(\cdot)\bigr)
 \nu(ds,dx),\\
S_{n} &\leq & t\leq T<S_{n+1}\quad\Rightarrow\quad
Y_t=\xi+\int_t^{T}\!\int_Ef\bigl(s,x,Y_s,Z_s(\cdot)\bigr) \nu(ds,dx).
\end{eqnarray*}
Using the form $Y\equiv(y^^n)$, and $Z\equiv(z^n)$ as defined by \eqref{B11},
this is equivalent to having for $\om\notin\n$
\begin{eqnarray}
&& S_{n}(\om) \leq   t<S_{n+1}(
\om)\leq T
\nonumber
\\
&&\qquad\Rightarrow\quad y^n_{D_n(\om)}(t) =
y^n_{D_n(\om)}\bigl(S_{n+1}(\om)\bigr) \nonumber
\nonumber
\\[-8pt]
\label{A1}
\\[-8pt]
\nonumber
&&\qquad\qquad\qquad\qquad\quad{}+\int
_t^{S_{n+1}(\om)}\!\int_Ef\bigl\{
\wy^{n+1}\bigr\}^n_{D_n(\om)}\bigl(s,x,y^n_{D_n(\om)}(s),
\nonumber\\
&&\qquad\qquad\qquad\quad\qquad{}-y^n_{D_n(\om)}(s)\bigr) \nu^n_{D_n(\om)}(ds,dx),
\nonumber\\
&& S_{n}(\om) \leq   t\leq
T<S_{n+1}(\om)
\nonumber
\\
&&\qquad\Rightarrow \quad y^n_{D_n(\om)}(t)=u^n_{D_n(\om)}\nonumber
\nonumber
\\[-8pt]
\label{A2}
\\[-8pt]
\nonumber
&&\quad\qquad\qquad\qquad\qquad
 {}+\int_t^{T}
\!\int_Ef\bigl\{\wy^{n+1}\bigr\}^n_{D_n(\om)}
\bigl(s,x,y^n_{D_n(\om)}(s), \\
&&\quad\qquad\qquad\qquad\qquad{}-y^n_{D_n(\om)}(s)
\bigr) \nu^n_{D_n(\om)}(ds,dx).\nonumber
\end{eqnarray}

Thus, if $Y$ is a solution and $\om\notin\n$, the function $y^n_{D_n(\om)}$
satisfies the differential equation in \eqref{A2} on the interval
$[S_n(\om)\wedge T,T]$, hence also on the interval $[0,T]$
because\vspace*{1pt} $\nu^n_{D_n(\om)}([0,S_n(\om)]\times E)=0$ and $y^n_{D_n(\om)}(t)
=y^n_{D_n(\om)}(S_n(\om))$ if $t\leq S_n(\om)$ and also $u^n_{D_n(\om)}=0$
and $y^n_{D_n(\om)}(t)=0$ if $S_n(\om)>T$: we thus have \eqref{B13}.

Conversely, assume that outside a null set $\n$ we have \eqref{B13} for
all $n$. Then obviously \eqref{A2} holds, and \eqref{A1} as well by
taking the difference\vspace*{1pt} $y^n_{D_n(\om)}(t)-y^n_{D_n(\om)}(S_{n+1}(\om))$.
Therefore, $Y$ solves the BSDE. This proves the first claim.

Assume further $\PP(S_{M+1}=\infty)=1$. Outside a null set, we have
$S_n=\infty$ for all $n>M$, so \eqref{B13} is trivially satisfied (with
both members equal to $0$) if $n>M$, and it reduces to \eqref{B14} when
$n=M$ because then $\nu^M_{D_M(\om)}([0,T]\times E)=0$, hence the second
claim.
\end{pf}

Equation \eqref{B13} leads us to consider the following equation with unknown
function $y$, for any given $n$,
\begin{eqnarray}
\label{B15}
&&\quad y(t)=u^n_{D}+\int_t^T
\!\int_Ef\{\wy\}^n_D
\bigl(s,x,y(s),-y(s)\bigr) \nu^n_{D}(ds,dx),\qquad t\in[0,T],\hspace*{-24pt}
\end{eqnarray}
where $D\in H_n$ is given, as well as the Borel function $\wy$ on
$[0,T]\times E$ with further $\wy(t,x)=0$ if $t\leq D^{\max}$. When
$D^{\max}=\infty$, and in view of our
prevailing convention $u^D=0$, plus $\nu^n_D([0,T]\times E)=0$ in this
case, this reduces to $y(t)=0$. Otherwise, this equation is a backward
ordinary integro-differential equation, and we have the following.

\begin{lem}\label{LB2} Equation \eqref{B15} has at most one solution,
and it has one as soon as
\begin{eqnarray}
\label{B16}
&&\int_0^T\!\int
_E\bigl|\wy(s,x)\bigr| \nu^n_D(ds,dx)<\infty.
\end{eqnarray}
In this case, the unique solution $y$ satisfies, for all $\rho\geq L+L'$,
\begin{eqnarray}
\bigl|y(t)\bigr| e^{\rho a^n_D(t)} &\leq &  \bigl|u^n_D\bigr|
e^{\rho a^n_D(T)}
\nonumber
\\[-3pt]
\label{B17}
\\[-13pt]
\nonumber
&&{}+ \int_t^T \!\int
_E \bigl(\bigl|f\{0\}^n_D(s,x,0,0)\bigr|+L\bigl|
\wy(s,x)\bigr| \bigr) e^{\rho a^n_D(s)} \nu^n_D(ds,dx)
\end{eqnarray}
and also, if $\rho>L+L'$ and with a constant $\BC$ depending only
on $(\rho,L,L')$,
\begin{eqnarray}
&& \int_t^T
\bigl|y(s)\bigr| e^{\rho a^n_D(s)} \,da^n_D(s)
\nonumber
\\
\label{B18}
&&\qquad\leq \BC \biggl(\bigl|u^n_D\bigr|e^{\rho a^n_D(T)}+\int
_t^T \!\int_E \bigl(\bigl|f\{0
\}^n_D(s,x,0,0)\bigr|\\
&&\qquad\quad {}+L\bigl|\wy(s,x)\bigr| \bigr) e^{\rho a^n_D(s)}
\nu^n_D(ds,dx) \biggr).\nonumber
\end{eqnarray}
\end{lem}

\begin{pf}
We\vspace*{1pt} have $f\{\wy\}^n_D(s,x,y(s),-y(s))=g(s,x,y(s))$, where
$g$ is a Borel function on $[0,T]\times E\times\R$, which by \eqref{B9}
satisfies
\begin{eqnarray*}
&& \bigl|g\bigl(s,x,y'\bigr)-g(s,x,y)\bigr|\leq \bigl(L+L'\bigr)\bigl|y'-y\bigr|,\\
&&
\int_0^T\!\int_E\bigl|g(s,x,0)\bigr| \nu^n_D(ds,dx)\\
&&\qquad\leq
\int_0^T\bigl|f\{0\}^n_D(t,x,0,0)\bigr| \nu^n_D(dt,dx)+
L\int_0^T\!\int_E\bigl|\wy(s,x)\bigr| \nu^n_D(ds,dx).
\end{eqnarray*}
The Lipschitz property of $g$ implies the uniqueness, and the existence
is classically implied by the finiteness of $\int_0^T\!\int_E|g(s,x,0)|
\nu^n_D(ds,dx)$, which holds under \eqref{B16} because of the last condition
in \eqref{B9}.

Next, under \eqref{B16}, the proof of the estimates is the same as in Lemma~\ref{LP2}. Namely, there is no jump here, so \eqref{P5} is replaced by
\begin{eqnarray*}
&& \bigl|y(t)\bigr|e^{\rho a^n_D(t)}+\rho\int_t^T\bigl|y(s)\bigr|e^{\rho a^n_D(s)}
\,da^n_D(s)\\
&&\qquad\leq\bigl|u^n_D\bigr|e^{\rho a^n_D(T)}
+\int_t^T\!\int_E \bigl(\bigl|g(s,x,0)\bigr|+\bigl(L+L'\bigr)\bigl|y(s)\bigr| \bigr)
e^{\rho a^n_D(s)} \nu^n_D(ds,dx).
\end{eqnarray*}
Note that here $\int_t^T|y(s)|e^{\rho a^n_D(s)} \,da^n_D(s)<\infty$ because
$a^n_D(T)<\infty$. We readily get~\eqref{B17} if $\rho\geq L+L'$, and \eqref{B18} if $\rho>L+L'$.
\end{pf}

We end this subsection with a technical lemma.

\begin{lem}\label{LB3}
For any $n\geq0$ and any nonnegative
Borel function $g$ on $[0,T]\times E\times H_n\times H_{n+1,}$ we have
\begin{eqnarray}
&& \int_0^T\!\int
_Eg\bigl(s,x,D_n,D_n\cup\bigl
\{(s,x)\bigr\}\bigr) \nu^n_{D_n}(ds,dx)
\nonumber
\\[-8pt]
\label{B19}
\\[-8pt]
\nonumber
&&\qquad = \E \bigl(g(S_{n+1},X_{n+1},D_n,D_{n+1})
e^{a^n_{D_n}(S_{n+1})} 1_{\{S_{n+1}\leq T\}}\vert \f_{S_n} \bigr).
\end{eqnarray}
Moreover,
the set $C'= \{D\in\h_n:~\int_0^T\!\int_E1_{\{D\cup\{(s,x)\}\in C\}}
 \nu^n_D(ds,dx)>0 \}$ is $\La_n$-negligible, if $C\subset H_{n+1}$ is $\La_{n+1}$-negligible.
\end{lem}

\begin{pf}
In view of \eqref{B7}, the left-hand side of \eqref{B19} is
\[
\int_0^T\!\int_Eg\bigl(s,x,D_n,D_n\cup\bigl\{(s,x)\bigr\}\bigr) e^{a^n_{D_n}(s)} G^n_{D_n}(ds,dx),
\]
so the first claim follows from the fact that $G^n_{D_n}$ is
the $\f_{S_n}$-conditional law of $(S_{n+1}X_{n+1})$. For the last claim, it
suffices to take the expectation of both sides of~\eqref{B17} with $g=
1_{[0,T]\times E\times H_n\times C}$: the right-hand side becomes
$\E (e^{a^n_{D_n}(S_{n+1})} \times 1_C(D_{n+1}) 1_{\{S_{n+1}\leq T\}} )$,
which vanishes because
$\La_{n+1}(C)=0$, whereas the left-hand side is positive if $\La_n(C')>0$.
\end{pf}

\textit{An example of an explicit solution}:  We will prove Theorem~\ref{TS1}
later, but here we show how Lemma~\ref{LB1} allows us to give an explicit solution,
in a special (but nontrivial) case of this theorem, with $M=2$.

We consider a state space $E=\{x_1,x_2,x_3\}$ with three elements and
suppose that $S_n=\infty$ for $n\geq3$ and that $X_1=x_1$
if $S_1<\infty$, whereas conditionally on $(S_1,S_2)$ and if $S_2<\infty$
then $X_2$ takes the two values $x_2$ and $x_3$ with probability $\frac{1}2$.
The law of the point process is thus completely characterized by the
law $H^1(dt)$ of $S_1$, and by the conditional law $H^2(S_1,dt)$ of $S_2$
knowing $S_1$ (so $H^2(s,dt)$ is a transition probability from $[0,\infty]$
into itself, satisfying $H^2(\infty,\{\infty\})=1$ and $H^2(s,(s,\infty])
=1$ if $s<\infty$). We also assume \ref{aa}, which amounts to the facts that
$H^1$ and $H^2_s$ have no atom except $\infty$, plus $H^1((T,\infty])>0$
and $H^2(s,(T,\infty])>0$.

We consider the linear equation
\begin{eqnarray}
&& Y_t+ \int_{(t,T]} \int
_E Z(s,x) \mu(ds,dx)
\nonumber
\\[-8pt]
\label{EE1}
\\[-8pt]
\nonumber
&&\qquad= 1_{\{S_2\leq T,X_2=x_2\}} +\int
_{(t,T]}\int_E Z(s,x) \nu(ds,dx).
\end{eqnarray}

With the notation \eqref{B2} and $\De=x_1$, say, we have $D_0=(0,x_1)$
and $D_1=((0,x_1),(S_1,x_1))$ reduces to $S_1$. Thus, we may take
\begin{eqnarray*}
u^0_D &=& 0,\qquad  u^1_D=0,\qquad
u^2_{D_2}=1_{\{S_2\leq T,X_2=x_2\}},\\
G'^0_D &=& H^1,\qquad  \phi^0_{D,t}=\ep_{x_1},\qquad  G'^1_{D_1}=H^2(S_1,\cdot),\qquad
  \phi^1_{D_1,t}=\tfrac{1}2 (\ep_{x_1}+\ep_{x_2}),\\
a^0_{D_0}(t) &=& a^0(t)=-\log H^1\bigl((t,\infty]\bigr), \\
a^1_{D_1}(t)&=&a^1_{S_1}(t)=
-\log H^2\bigl(S_1,(t,\infty]\bigr).
\end{eqnarray*}
Moreover, in \eqref{B5} $y^0_D(t)$ is a function $y^0(t)$, and
$y^1_{D_1}(t)$ takes the form $y^1_{S_1}(t)\times 1_{\{S_1\leq T\}}$ for some\vspace*{1pt}
function\vspace*{1pt} $(r,t)\mapsto y_r^1(t)$ on $[0,T]^2$, whereas by \eqref{B14} we
may take $y^2_{D_2}(t)=u^2_{D_2}$ for all $t$. The form\vspace*{1pt} of the generator
implies that in \eqref{B8} we have $f\{z^n\}^n_{D_n}(t,x,y,w)=
w-z^n_{D_n}(t,x)$. Then,\vspace*{1pt} writing \eqref{B13} for $n=1$ and $n=0$ gives us
(below, $r$ stands for $S_1$)
\begin{eqnarray*}
y^1_r(t) &=& \frac{1}2  \bigl(a^1_r(T)-a^1_r(t) \bigr)-\int_t^T y^1_r(s) \,da^1_r(s),\\
y^0(t) &=& \int_t^Ty^1_s(s) \,da^0(s)-\int_t^T y^0(s) \,da^0(s).
\end{eqnarray*}
This is a system of linear ODEs, whose explicit solution is [recall
$a^1_s(s)=0$]
\begin{eqnarray*}
y^1_r(t) &=& \tfrac{1}2  \bigl(1-e^{a^1_r(t)-a^1_r(T)} \bigr),\\
y^0(t)&=& \frac{1}2\int_t^Te^{a^0(t)-a^0(s)}  \bigl(1-e^{-a^1_s(T)} \bigr)
 \,da^0(s).
\end{eqnarray*}
Upon\vspace*{1.5pt} replacing $a^1_s(t)$ and $a^0(t)$ by $-\log \BH^2_s(t)$ and
$-\log \BH^1(t)$, and using $y^2_{D_2}(t)=1_{\{S_2\leq T,X_2=x_2\}}$, we
obtain the following explicit form for the unique solution:
\begin{eqnarray*}
&& t\in[S_2,T] \quad\Rightarrow\quad Y_t=1_{\{X_2=x_2\}},\\
&& t\in[S_1\wedge T,S_2\wedge T]\quad\Rightarrow\quad Y_t=\frac{H^2(S_1,(t,T])}{
2\BH^2(S_1,(t,\infty])},\\
&& t\in[0,S_1\wedge T]\quad\Rightarrow\quad Y_t=\frac{1}{2H^1((t,\infty])}\int_t^T
H^2\bigl(s,(s,T]\bigr) H^1(ds).
\end{eqnarray*}

\begin{rema}
In this example, we have \eqref{S20} with $M=2$, so the
uniqueness holds by Theorem~\ref{TS1}. We also have \eqref{S21}, but the
process $A$ is not necessarily bounded: nevertheless we do have
existence.
\end{rema}

\subsection{Some counter-examples when \texorpdfstring{\protect\ref{aa}}{(A)} fails}\label{ss-3-3}

In all the paper, we assume \ref{aa}, and it is enlightening to see what happens
when this assumption fails. We are not going to do any deep study of
this case, and will content ourselves with the simple situation where the
point process is univariate and has a single point, that is, $E=\{\De\}$ is
a singleton, and
\[
N_t = 1_{\{S\leq t\}},
\]
where $S$ is a variable with values in $(0,T]\cup\{\infty\}$. The filtration
$(\f_t)$ is still the one generated by $N$, and $G$ denotes the law of $S$,
whereas $g(t)=G(t,\infty]$: those are the same as in \eqref{B6}, in our
simplified setting.

The equation is \eqref{S4}, but since $A_t=A_{t\wedge S}$ and any predictable
process is nonrandom, up to time $S$, it now reads as
\begin{eqnarray}
\label{B20}
&&
Y_t+Z_S 1_{\{t<S\leq T\}}=\xi+\int
_{(t,S\wedge T]}f(s,Y_{s-},Z_s) \,dA_s,
\end{eqnarray}
with $f$ a Borel function on $[0,T]\times\R\times\R$, Lipschitz in its
last two arguments, and  such that $\int_0^T|f(s,0,0)| \,dA_s<\infty$.

Assumption \ref{aa} fails if (A$_1$) or (A$_2$) or both fail.
Below, we examine what happens if either one of these two partial assumptions
fails.
\begin{longlist}[(2)]
\item[(1)]  \textit{When $G$ has an atom}.  Here, we assume that (A$_1$) does not hold,
that is, $A$ is discontinuous, whereas $\PP(S=\infty)>0$, so (A$_2$) holds. We
will see that in this case the existence of a solution to \eqref{B20} is
not guaranteed.

To see this, we consider the special case where $S$ only takes the two values
$r\in(0,T]$ and $\infty$, with respective positive probabilities $p$
and $1-p$. We have $N_t=1_{\{r\leq t\}} 1_{\{S=r\}}$
and $A_t=p1_{\{t\geq r\}}$, so
only the values of $f(t,y,z)$ at time $t=r$ are relevant, and we may assume that
$f=f(y,z)$ only depends on $y,z$. Note also that $\xi$ takes the form
\[
\xi=a 1_{\{S=r\}}+b 1_{\{S=\infty\}}\qquad \mbox{where $a,b\in\R$}.
\]
Moreover, only the value $Z_r(\om)$ is involved, and it is nonrandom,
and any solution $Y$ is constant on $[0,r)$ and on $[r,T]$, that is,
we have for $t\in[0,T]$
\begin{eqnarray}
&& Z_r(\om)=\ga, \qquad Y_t=\de 1_{\{t<r\}}+\rho 1_{\{t\geq r,
S=r\}}+\eta 1_{\{t\geq r, S=\infty\}}\nonumber\\
\eqntext{\mbox{where $\ga,\de,\rho,\eta\in\R$}.}
\end{eqnarray}

Here, $a,b$ are given, and $\ga,\de,\rho,\eta$ constitute the ``solution'' of
\eqref{B20}, which reduces to the four equalities
\[
\eta=b,\qquad  \rho=a, \qquad \de=b+pf(\de,\ga),
\qquad \de+\ga=a+pf(\de,\ga),
\]
which in turn give us
\[
\ga=a-b,\qquad  \de=b+pf(\de,a-b).
\]
The problem is that the last equation may not have a solution, and if it
has one it is not necessarily unique. For example, we have:
\begin{eqnarray*}
&&\mbox{if $f(y,z)=\dfrac{1}p \bigl(y+g(z)\bigr)$},\\
&&\qquad\mbox{then }
\cases{\mbox{if $a+g(a-b)=0$ there are infinitely many solutions},\vspace*{3pt}\cr
\mbox{if $a+g(a-b)\neq0$ there is no solution.}
}
\end{eqnarray*}
\item[(2)] \textit{When $G$ is supported by $[0,T]$}.  Here, we suppose that $G$ has
no atom, but is supported by $[0,T]$. This corresponds to having (A$_1$),
but not (A$_2$), and we have $A_t=a(t\wedge S)$, where
$a(t)=-\log g(t)$ is increasing, finite for $t<v$ and infinite if
$t\geq v$, where $v=\inf(t: g(t)=0)\leq T$ is the right end point of the
support of the measure $G$.

We will also consider a special
generator, and more specifically the equation
\begin{eqnarray}
\label{B21}&& Y_t+\int_{(t,T]}Z_s
(dN_s-dA_s)=\xi,
\end{eqnarray}
which is \eqref{S6} with $f\equiv0$, and \eqref{S4} with $f(t,y,z)=z$.

When $\xi$ is integrable, the martingale representation theorem for point
processes yields that $\xi=\E(\xi)+\int_0^TZ_s(dN_s-dA_s)$ for some
predictable and $dA_t$-integrable process $Z$, hence $Y_t=\E(\xi\vert \f_t)$
is a solution. But this is \textit{not} the only one. Indeed, recalling that
here $\xi=h(S)$ is a (Borel) function of $S$, we have
the following.
\end{longlist}

\begin{prop}\label{PB1}
Assume that $\PP(S\leq T)=1$ and that the law of
$S$ has no atom, and also that $\xi$ is integrable. Then a process $Y$ is
a solution of \eqref{B21} if and only if, outside a $\PP$-null set, it takes
the form
\begin{eqnarray}
\label{B22}
&& Y_t=\xi 1_{\{t\geq S\}}+ \biggl(w-\int
_0^t e^{-A_s} h(s) \,dA_s
\biggr) e^{A_t} 1_{\{t<S\}}
\end{eqnarray}
for an arbitrary real number $w$, and the associated process $Z$ can be taken
as $Z_t=h(t)-Y_{t-}$.
\end{prop}

Note that $Y_0=w$ in \eqref{B22}, so in particular it follows that \eqref{B21}
has a unique solution for any initial condition $Y_0=w\in\R$. This is in
deep contrast with Theorems \ref{TS1} or \ref{TS2}, and it holds even
for the trivial case $\xi\equiv0$: in this trivial case, $Y_t=0$ is of
course a solution, but $Y_t=we^{A_t} 1_{\{t<S\}}$ for any $w\in\R$ is also
a solution.

\begin{pf*}{Proof of Proposition~\protect\ref{PB1}}
Any solution $(Y,Z)$ satisfies $Y_t=\xi$ if $t\geq S$
and $Y_t=y(t)$ if $t<S$, where $y$ is a continuous (nonrandom) function
on $[0,v)$ (recall that $S<v$ a.s., and ess sup $S=v$). Since further
\eqref{S16} holds, one may always take the associated predictable process
$Z$ to be $Z_t=h(t)-Y_{t-}$. Then writing \eqref{B21} for $t=0$ and $t$
arbitrary in $[0,v)$, we see that $Y$ is a solution if and only if
\[
t\in[0,v)\quad\Rightarrow\quad y(t)=y(0)+\int_0^t\bigl(y(s)-h(s)\bigr) \,da(s).
\]
This is a linear ODE  whose solutions are exactly the functions
\[
y(t)= \biggl(w-\int_0^te^{-a(s)} h(s) \,da(s) \biggr) e^{a(t)}
\]
for $w\in\R$ arbitrary
[since $\int_0^t|h(s)| \,da(s)\leq \frac{1}{g(t)} \E(|\xi|)$ is finite for all
$t<v$]. This completes the proof.
\end{pf*}

\begin{rema}
The previous result does not depend on the special form of the generator $f$,
in the sense that for any $f$ satisfying \eqref{S5} and under the assumptions
of Proposition~\ref{PB1}, for any $w\in\R$ the BSDE admits a unique solution
starting at $Y_0=w$: of course an explicit form such as \eqref{B22} is no
longer available, but the proof of this result follows exactly the same
argument as above.
\end{rema}

\begin{rema}
Jeanblanc and R\'eveillac \cite{JR14} have studied some
cases of BSDEs driven by a Wiener process, for which the generator
``explodes'' at the terminal time $T$. This bears some resemblance with
the previous setting, in which $A_t=a(t\wedge S)$ and $a(t)\to\infty$ as
$t\to T$. They show for example that, in the affine case, and under
appropriate assumptions, there is no solution when $\PP(\xi\neq0)>0$,
and infinitely many solutions when $\xi\equiv0$. Of course, the setting
is quite different (a Wiener process instead of a point process), so the
results are not really comparable, but they find cases like when (A$_2$)
fails (no solutions) and like when (A$_1$) fails (infinitely many
solutions).
\end{rema}

\section{Proof of the main results}\label{s-4}

We start with an auxiliary lemma  needed for proving the existence of a
solution.

\begin{lem}\label{LC1} Assume \eqref{S20} and that $A_T\leq K$ for some
constant $K$. Let $m\in\{1,\ldots,M\}$, and suppose that we have
$y^n_{D_n}(t)$ for $n=m,m+1,\ldots,M$, such that~\eqref{B13} holds if
$m\leq n<M$ and \eqref{B14} holds if $n=M$, outside a null set. Then for
$n$ between $m$ and $M-1$, we have the (rather coarse) estimate
\begin{eqnarray}
v_n &:=& \int_0^T
\!\int_E\bigl(\bigl|f\{0\}^n_{D_n}(s,x,0,0)\bigr|+L\bigl|y^{n+1}_{D_n\cup\{(s,x)\}}(s)\bigr|
\bigr) \nu^n_{D_n}(ds,dx)
\nonumber
\\
\label{C1}
&\leq &  (1+L)^M e^{MK(2+L+L')} \\
\nonumber
&&{}\times\E \biggl(\int_{S_n\wedge T}^T
\int_E \bigl|f(s,x,0,0)\bigr| \nu(ds,dx)+|\xi| 1_{\{S_n\leq T\}}\vert
\f_{S_n} \biggr).
\end{eqnarray}
\end{lem}

\begin{pf}
(1) We first prove that $A_T\leq K$ implies
\begin{eqnarray}
\label{C2}
&& n\geq 0,\qquad D\in H_n\quad\Rightarrow\quad a^n_{D}(T)
\leq K
\end{eqnarray}
for a suitable version of the $a^n_D$'s, which amounts to proving
$a^n_{D_n}(T)\leq K$ a.s. To check this, we observe that for any $\ga>1$
\begin{eqnarray*}
e^{(\ga-1)a^n_{D_n}(T)}&=& \E \bigl(e^{\ga a^n_{D_n}(T)} 1_{\{S_{n+1}>T\}}
\vert \f_{S_n} \bigr)\\
&=&\E \bigl(e^{\ga a^n_{D_n}(T\wedge S_{n+1})} 1_{\{S_{n+1}>T\}}
\vert \f_{S_n} \bigr)\leq e^{K\ga},
\end{eqnarray*}
because $a^n_{D_n}(T\wedge S_{n+1})\leq A_T$ by \eqref{B7}. This
implies $a^n_{D_n}(T)\leq \frac{K\ga}{\ga-1}$ a.s. and, being true for
all $\ga>1$, it yields \eqref{C2}.

(2) By Lemma~\ref{LB3}, we have outside a null set
\[
v_n=\E \bigl(e^{a^n_{D_n}(S_{n+1})} \bigl(\bigl|f\{0\}^n_{D_n}(S_{n+1},X_{n+1},0,0)\bigr|
+L\bigl|y^{n+1}_{D_{n+1}}(S_{n+1})\bigr|\bigr) 1_{\{S_{n+1}\leq T\}}\vert \f_{S_n} \bigr).
\]
Equation \eqref{B8} yields $f\{0\}^n_{D_n}(t,x,0,0)=f(t,x,0,0)$ if
$S_n<t\leq S_{n+1}$, whereas $u^n_{D_n}=0$ if $S_n>T$, and $u^n_{D_n}=\xi$ if
$S_n\leq T<S_{n+1}$. In view of \eqref{B14} and \eqref{C2}, we first deduce
\[
v_{M-1}\leq e^K \E \bigl(\bigl(\bigl|f(S_M,X_M,0,0)\bigr|+L|\xi|
\bigr) 1_{\{S_M\leq T\}}\vert \f_{S_{M-1}} \bigr).
\]
It also gives us for $n\leq M-2$, upon using \eqref{B17} with
$n+1$ and $\rho=L+L'$, and~\eqref{C2} again
\begin{eqnarray*}
v_n &\leq& e^K \E \bigl( \bigl(\bigl|f(S_{n+1},X_{n+1},0,0)\bigr|
+Le^{(L+L')K}\bigl(\bigl|u^{n+1}_{D_{n+1}}\bigr|
+v_{n+1}\bigr) \bigr) 1_{\{S_{n+1}\leq T\}}
\vert \f_{S_n} \bigr)\\
&\leq & e^K \E \bigl( \bigl(\bigl|f(S_{n+1},X_{n+1},0,0)\bigr|
\\
&&{}+Le^{(1+L+L')K}\bigl(|\xi| 1_{\{S_{n+2}>T\}}
+v_{n+1}\bigr) \bigr) 1_{\{S_{n+1}\leq T\}}
\vert \f_{S_n} \bigr),
\end{eqnarray*}
where we have used $\PP(S_{n+2}>T\vert \f_{S_{n+1}})\geq e^{-K}$, which implies
\begin{eqnarray*}
\E \bigl(|\xi| 1_{\{S_{n+2}> T\ge S_{n+1}\}}\vert \f_{S_{n+1}} \bigr)
&=& \E \bigl(\bigl|u^{n+1}_{D_{n+1}}\bigr|
 1_{\{S_{n+2}> T\ge S_{n+1}\}}\vert \f_{S_{n+1}} \bigr)\\
&=&\bigl|u^{n+1}_{D_{n+1}}\bigr| 1_{\{ T\ge S_{n+1}\}} \PP(S_{n+2}>T\vert \f_{S_{n+1}})
\\
&\geq& \bigl|u^{n+1}_{D_{n+1}}\bigr| 1_{\{ T\ge S_{n+1}\}} e^{-K}.
\end{eqnarray*}
Iterating the estimates for $v_n$, and by successive conditioning, we deduce
\begin{eqnarray*}
v_n &\leq &  (1+L)^M  e^{MK(2+L+L')}\\
&&{}\times
      \E \biggl(\sum_{i=n}^{M-1}
\bigl|f(S_{i+1},X_{i+1},0,0)\bigr| 1_{\{S_{i+1}\leq T\}}+L|\xi|
 1_{\{S_i\leq T<S_{i+1}\}}\vert \f_{S_n} \biggr)\\
 & \leq &  (1+L)^M  e^{MK(2+L+L')}\\
 &&{}\times \E \biggl(\int_{S_n\wedge T}^T\int_E
\bigl|f(s,x,0,0)\bigr| \mu(ds,dx)+L|\xi| 1_{\{S_n\leq T\}}\vert \f_{S_n} \biggr).
\end{eqnarray*}
Since $\nu$ is the compensator of $\mu$, this is equal to the right-hand
side of \eqref{C1}, hence the result.
\end{pf}

\begin{pf*}{Proof of Theorem~\protect\ref{TS1}}
(a) We first prove the uniqueness. Let
$Y\equiv(y^n)$ and $Y'\equiv(y'^n)$ be two solutions.
By Lemma~\ref{LB1}, for any $n=0,\ldots,M$ we have a subset $B_n$ of
$H_n$ with $\La_n(B_n^c)=0$ and such that for any $D\in B_n$ both $y^n_D$
and $y'^n_D$ satisfy~\eqref{B13} if $n<M$ and \eqref{B14} if $n=M$.

The proof is done by downward induction. The induction hypothesis $K(n)$
is that for all $m=n,\ldots,M$ we have a subset $B(n,m)$ of $H_m$ with
$\La_m(B(n,m)^c)=0$ such that $y^m_D\equiv y'^m_D$ for all $D\in B(n,m)$.
That $K(M)$ holds with $B(M,M)=B_M$ is obvious, and $K(0)$
yields $Y_t=Y'_t$ a.s. for all $t$.

It remains to show that $K(n+1)$ for some $n$ between $0$ and $M-1$
implies $K(n)$. Assuming $K(n+1)$, we set $B(n,m)=B(n+1,m)$ for $m>n$ and
let $B(n,n)$ be the intersection of $B_n$ and of the set of all $D\in H_n$
such that $y^{n+1}_{D\cup\{(s,x)\}}=y'^{n+1}_{D\cup\{(s,x)\}}$ for
$\nu^n_D$-almost all $(s,x)$. By virtue of the last claim in Lemma~\ref{LB3}
applied with $C=B(n+1,n+1)^c$, plus $\La_n(B_n^c)=0$, we have
$\La_n(B(n,n)^c)=0$. Then Lemma~\ref{LB2} yields
$y^n_D=y'^n_D$ when $D\in B(n,n)$, hence $K(n)$ holds.

(b) We now turn to the existence, assuming further $A_T\leq K$ and \eqref{S21}.
We construct the family $(y^n_{D_n}(t))$ by downward induction on $n$, starting
with $y^M_D(t)=u^M_D$ for all $D\in H_M$, hence \eqref{B14} holds everywhere.
Suppose now that we have a null set $C_{n+1}$ and functions $y^m_{D_m}$ for
$m=n+1,\ldots,M-1$, each one satisfying~\eqref{B13} outside $C_{n+1}$.
The assumption \eqref{S21} and Lemma~\ref{LC1} imply $\E(v_n)<\infty$,
so the set $C_n=C_{n+1}\cup\{v_n=\infty\}$ is negligible. Now,
\eqref{B13} is the same as \eqref{B15}
with $D=D_n$ and $\wy(s,x)=y^{n+1}_{D_n\cup\{(s,x)\}}(s) 1_{\{t>D^{max}\}}$,
which is well
defined for $G^n_{D_n}$-almost all $(s,x)$, hence for $\nu^n_{D_n}$-almost
all $(s,x)$. Therefore, outside $C_n$ these terms satisfy \eqref{B16}, and it
follows that \eqref{B15} has a unique solution $y^n_{D_n}$. This validates
the induction, hence \eqref{S11} has a solution, necessarily a.s. unique
by part (a) above.

(c) It remains to prove the last claims. We denote by $Y$ the (a.s. unique)
solution, and recall that the associated predictable function $Z$ can be
chosen as $Z\equiv(z^n)$ with the form \eqref{B1}. Since $N_T\leq M$, the
last two claims amount to proving that $\E(U_n)<\infty$ for all $n\leq M$,
where $U_n=\int_{S_n\wedge T}^{S_{n+1}\wedge T}\int_E(|Y_s|+|Z(s,x)|)
\nu(ds,dx)$. Since $U_M=0$ because $A_T=A_{T\wedge S_M}$, we restrict our
attention to the case $n<M$. \eqref{B4}, \eqref{B7} and \eqref{B1} yield
$U_n\leq2V_n+W_n$, where
\begin{eqnarray*}
V_n &=& \int_{S_n\wedge T}^T \bigl|y^n_{D_n}(s)\bigr| \,da^n_{D_n}(s),\qquad
W_n=\int_{S_n\wedge T}^T\int_E \bigl|y^{n+1}_{D_n\cup\{(s,x)\}}(s)\bigr|
\nu^n_{D_n}(ds,dx).
\end{eqnarray*}
On the one hand, $LW_n\leq v_n$, so \eqref{S21}
and \eqref{C1} yield $\E(W_n)<\infty$. On the other hand, applying\vspace*{1.5pt}
first \eqref{B18} with any $\rho>L+L'$ and \eqref{C2} and then
$\PP(S_{n+1}>T\vert \f_{S_{n}})\geq e^{-K}$ and \eqref{C1}, we get
\begin{eqnarray*}
\E(V_n)&\leq & \BC e^{K(L+L')} \E \bigl(\bigl|u^n_{D_n}\bigr|
1_{\{S_n\leq T\}}+v_n \bigr)
\\
&\leq & \BC e^{K(1+L+L')} \E \bigl(|\xi|
1_{\{S_n\leq T<S_{n+1}\}}+v_n \bigr)<\infty.
\end{eqnarray*}
This completes the proof.
\end{pf*}

\begin{pf*}{Proof of Theorem~\protect\ref{TS2}}
(a) The uniqueness has been proved in
Lemma~\ref{LP3}. For the existence, we will ``localize'' the problem in the
following way: for any $n\geq1$ we set $T_n=S_n\wedge\inf(t: A_t\geq n)$ and
we consider the equation
\begin{eqnarray}
 &&  Y^{(n)}_t+\int
_t^T\!\int_EZ^{(n)}(s,x)
\mu^{(n)}(ds,dx) \nonumber\\
&&\qquad=\xi^{(n)}+\int_t^T
\!\int_E f\bigl(s,x,Y_{s}^{(n)},Z^{(n)}_s(
\cdot)\bigr) \nu^{(n)}(ds,dx),
\nonumber
\\[-8pt]
\label{C3}
\\[-8pt]
\nonumber
&& \mu^{(n)}(ds,dx)=\mu(ds,dx) 1_{\{s\leq T_n\}}, \qquad\nu^{(n)}(ds,dx)=
\nu(ds,dx) 1_{\{s\leq T_n\}},\hspace*{-12pt}
\\
&&\xi^{(n)}=\xi 1_{\{T<T_n\}}.\nonumber
\end{eqnarray}
Then $\nu^{(n)}$ is the compensator of $\mu^{(n)}$, relative to $(\f_t)$ and
also to the smaller filtration $(\f^{(n)}_t=\f_{t\wedge T_n})$ generated by
$\mu^{(n)}$, whereas $\xi^{(n)}$ is $\f^{(n)}_T$-measurable. The two
marginal processes $N^{(n)}_t=\mu^{(n)}([0,t]\times E)$ and $A^{(n)}_t=
\nu^{(n)}([0,t]\times E)$ satisfy $A^{(n)}_T\leq n$ and $N^{(n)}_T\leq n$,
and \eqref{S23} clearly implies \eqref{S21} for $\xi^{(n)}$ and $\nu^{(n)}$.
Therefore, Theorem~\ref{TS1} implies the existence of an a.s. unique solution
$(Y^{(n)},Z^{(n)})$ to~\eqref{C3}, and the last claim of this theorem
further implies that $\|(Y^{(n)},Z^{(n)})\|_{\al,\be}^{(n)}<\infty$,
where the previous norm is the same as \eqref{S22} with $(A,N,\nu)$
substituted with $(A^{(n)},N^{(n)},\nu^{(n)})$.

For $n'>n$, set
\begin{eqnarray*}
\BY^{(n,n')}&=& \sup_{s\in[0,T]}
 \bigl(e^{\be A_{s}} \al^{N_{s}}\bigl|Y^{(n')}_s-Y^{(n)}_s\bigr| \bigr),\\
 \w^{(n,n')}_{(s,t]} &=&\w^{\al,\be}_{(s,t]}\bigl(Y^{(n')}-Y^{(n)},
Z^{(n')}-Z^{(n)}\bigr),
\end{eqnarray*}
the latter being computed as in \eqref{P2} with $(A,N,\nu)$.

We now proceed to bound these variables, and to this end we observe that
\begin{eqnarray*}
&& Y^{(n')}_{T_n\wedge t}+\int_t^T\!\int_EZ^{(n')}(s,x) \mu^{(n)}(ds,dx)
\\
&&\qquad=Y^{(n')}_{T_n\wedge T}+\int_t^T\!\int_Ef\bigl(s,x,Y_{T_n\wedge s}^{n'},
Z^{(n')}_s(\cdot)\bigr) \nu^{(n)}(ds,dx),
\end{eqnarray*}
so $(Y^{(n')}_{T_n\wedge t},Z^{(n')})$ is a solution of \eqref{C3} with
terminal value $Y^{(n')}_{T_n\wedge T}$ instead of $\xi^{(n)}$, and clearly
has a finite $\|\cdot\|_{\al,\be}^{(n')}$ norm. It then follows from \eqref{P6}
and \eqref{P7}, plus the maximal inequality for
martingales, that for any $\ep>0$ we have
\begin{eqnarray}
&& \PP \Bigl(\sup_{t\in[0,T]}
e^{\be A_{T_n\wedge t}} \al^{N_{T_n\wedge t}}
 \bigl|Y^{(n')}_{T_n\wedge t}-Y^{(n)}_t\bigr|>\ep \Bigr)
\leq \frac{\de(n,n')}{\ep},\nonumber\\
&& \E\bigl(\w^{(n,n')}_{(0,T_n\wedge T]}\bigr)\leq C\de\bigl(n,n'\bigr)\nonumber\\
\eqntext{\displaystyle\mbox{where }\de\bigl(n,n'\bigr)=\E \bigl(\bigl|Y^{(n')}_{T_n\wedge T}-\xi^{(n)}\bigr|
e^{\be A_{T_n\wedge T}} \al^{N_{T_n\wedge T}} \bigr).}
\end{eqnarray}
If $T_n>T$, we have $Y^{(n')}_{T_n\wedge T}=Y^{(n')}_T=\xi^{(n')}=\xi
=\xi^{(n)}$, and otherwise $\xi^{(n)}=0$. Hence,  \eqref{P3} yields
\begin{eqnarray}
\de\bigl(n,n'\bigr) &=& \E \bigl(\bigl|Y^{(n')}_{T_n}\bigr|
e^{\be A_{T_n}} \al^{N_{T_n}} 1_{\{T_n\leq T\}} \bigr)
\leq \de_n\nonumber\\
\eqntext{\displaystyle \mbox{where }
\de_n=\E \biggl(|\xi| e^{\be A_T} \al^{N_T} 1_{\{T_n\leq T\}}+
\int_{T_n\wedge T}^T\!\int_E
\bigl|f(s,x,0,0)\bigr| e^{\be A_s} \al^{N_s} \nu(ds,dx)\biggr).}
\end{eqnarray}
If $T_n<t\leq T$, we have $Y^{(n)}_t=\xi^{(n)}=0$ and we may take
$Z^{(n)}(t,x)=0$, whereas if $T_{n'}\leq t\leq T$ we have
$Y^{(n')}_t=\xi^{(n')}=0$ and we may take $Z^{(n')}(t,x)=0$, hence
\begin{eqnarray*}
&& \w^{n,n'}_{(0,T]}-\w^{n,n'}_{(0,T_n\wedge T]}\\
&&\qquad=
\int_{T_n\wedge T}^{T_{n'}\wedge T}\int_E
 \bigl(\bigl|Y^{(n')}_s\bigr|+\bigl|Z^{(n')}(s,x)\bigr| \bigr)
e^{\be A_{s\wedge T_{n'}}} \al^{N_{s\wedge T_{n'}}} \nu^{(n')}(ds,dx).
\end{eqnarray*}
This and \eqref{P4} yield $\E (\w^{n,n'}_{(0,T]}-\w^{n,n'}_{(0,
T_n\wedge T]} )\leq C\de_n$. Gathering all those partial results, we
end up with
\begin{eqnarray}
\label{C4}
&& \PP\bigl(\BY^{(n,n')}>\ep\bigr)\leq \PP(T_n\leq
T)+\frac{\de_n}{\ep},\qquad \E\bigl(\w^{n,n'}_{(0,T]}\bigr)\leq 2C
\de_n.
\end{eqnarray}
In view of \eqref{S23} and the property $T_n\to\infty$ as $n\to\infty$, the
dominated convergence theorem implies $\de_n\to0$, hence both left sides
in \eqref{C4} go to $0$ as $n\to\infty$, uniformly in $n'>n$.
It follows that the sequence $Y^{(n)}$ is Cauchy for the convergence in
probability, in the Skorokhod space $\D([0,T])$ endowed with the uniform
metric, and that the pair $(Y^{(n)},Z^{(n)})$ is Cauchy in the space
$\laa_{\al,\be}^1$. Therefore, these sequences converge
in these spaces, to two limits $Y$ and $(Y',Z)$, with $Y$ c\`adl\`ag adapted
and $(Y',Z)\in\laa^1_{\al,\be}$ and $Z$ predictable and satisfying
$\int_0^T\int_E|Z(s,x)| \nu(ds,dx)<\infty$ a.s.; we can of course find versions
of the two limits for which $Y'=Y$ is the same process. Note that, since
all $Y^n$ are continuous outside the points $S_n$'s, the same is true
of $Y$.

We further deduce
$\E (\int_0^T\!\int_E|Z^{(n)}(s,x) -Z(s,x)| \nu(ds,dx) )
\to0$, implying $\E (\int_0^T\!\int_E|Z^{(n)}(s,x)-Z(s,x)| \mu(ds,dx) )
\to0$, and thus $\int_t^T\!\int_EZ^{(n)}(s,x) \mu(ds,\break dx)\toop
\int_t^T\!\int_EZ(s,x) \mu(ds, dx)$. Similarly, we obtain
$\int_t^T\!\int_Ef(s,x,Y^{(n)}_s,\break Z^{(n)}_s(\cdot))\nu(ds,dx)\toop\int_t^T\!\int_E
f(s,x,Y_s, Z_s(\cdot))\nu(ds,dx)$ (we use the Lipschitz property of
$f$ here), and of course $Y^{(n)}_t\toop Y_t$ for each $t$.
Since $(Y^{(n)},Z^{(n)})$ solves \eqref{C2}, by passing to the limit we
deduce that $(Y,Z)$ solves \eqref{S11}, and it clearly belongs to
$\laa^1_{\al,\be}$, thus ending the proof of the claim (a).

(b) We only need to prove that \eqref{S24} for some $\ep>0$
implies \eqref{S23} for all $\al>0$ and $\be\geq0$, when $A_T\leq K$ for some
constant\vspace*{1pt} $K$. Since $\nu([0,T]\times E)=A_T$ and $\al^{N_t}\leq
(\al\vee1)^{N_T}$ and $e^{\be A_t}\leq e^{\be K}$, by H\"older's inequality
it is clearly enough to show that $\al^{N_T}$ is in all $\LLL^p$ when
$\al>1$, or equivalently that $\E(\al^{N_T})<\infty$ for all $\al>1$.

We consider the nonnegative increasing process $U_t=\al^{N_t}$, which
satisfies the equation
\[
U_t=1+\al\int_0^tU_{s-} \,dN_s=1+\al\int_0^tU_{s-} \,dA_s+
\al\int_0^tU_{s-} (dN_s-dA_s).
\]
The last term is a local martingale, and a bounded martingale if we stop it
at time $S_n\wedge T$, because $N_{S_n}\leq n$ and $A_T\leq K$ and
$U_{t-}\leq\al^{n-1}$ if $t\leq S_n\wedge T$. Therefore, for any stopping
time $S\leq S'_n:=S_n\wedge T$ we have
\[
\E(U_{S-})\leq \E(U_S)=1+\al \E \biggl(\int_0^SU_{s-} \,dA_s \biggr).
\]
Then one applies the Gronwall-type lemma (3.39) in \cite{JM81} and
$A_{S'_n}\leq K$ to obtain that $\E(U_{S'_n-})\leq K'$ for a constant $K'$
which only depends on $K$ and $\al$. Letting $n\to\infty$ and using the
fact that $U_T\leq \al U_{T-}$, the monotone convergence theorem yields
$\E(U_T)\leq\al K'$ as well, hence the result.
\end{pf*}

\section{Application to a control problem}\label{s-5}

In this section, we show how what precedes can be put in use for solving
a control problem. As before, we are given the multivariate point process
$\mu$ of \eqref{S1} on $(\Om,\f,\PP)$, generating the filtration $(\f_t)$,
and satisfying \ref{aa}. The control problem is specified by the following data:
\begin{itemize}
\item a \textit{terminal cost}, which is an $\f_T$-measurable random
variable $\xi$;
\item an \textit{action} (or, decision) \textit{space}, which is a measurable
space $(U,\ua)$, and an associated predictable function $r$ on
$\Om\times[0,T]\times E\times U$, which specifies how the control acts;
\item a \textit{running cost}, which is a predictable function $l$ on
$\Om\times[0,T]\times U$.
\end{itemize}
These data should satisfy the following.
\renewcommand{\theass}{(B)}
\begin{ass}\label{bb}
There is a constant $C>0$ such that,
with $A$ and $N$ as in \eqref{S2} and \eqref{S3},
\begin{eqnarray}
\label{CC1}
&& 0\leq r(\om,t,x,u)\leq C,
\\
\label{CC2}
&& \E \bigl(e^{A_T} C^{N_T} \bigr)<\infty.
\end{eqnarray}
We also have, for two constants $\al\in[1,\infty)\cap(C,\infty)$ and
$\be>1+C$,
\begin{eqnarray}
&& \E \biggl(e^{\be A_T} \al^{N_T} |\xi|+\int
_0^Te^{\be A_s} \al^{N_s}\Bigl|\inf
_{u\in U} l(s,u)\Bigr| \,dA_s
\nonumber
\\[-8pt]
\label{CC3}
\\[-8pt]
\nonumber
&&\qquad{}+\int_0^Te^{A_s}
C^{N_s} \sup_{u\in U} \bigl|l(s,u)\bigr| \,dA_s
\biggr) < \infty.
\end{eqnarray}
\end{ass}

We denote by $\aaa$ the set of $U$-valued predictable processes. An element
of $\aaa$ is called an \textit{admissible  control}, and it operates as follows.
With $u=(u_t)\in\aaa$ we associate the probability measure $\PP_u$ on
$(\Omega,\f)$ which is absolutely continuous with respect to $\PP$ and
admits the density process
\begin{eqnarray}
&& L^u_t=\exp \biggl(\int_0^t\!\int_E \bigl(1-r(s,x,u_s)\bigr) \nu(ds,dx) \biggr)
\prod_{n\ge1 : S_n\le t}r({S_n},X_n,u_{S_n}),\nonumber
  \\
  \eqntext{t\in [0,T],}
  \end{eqnarray}
with the convention that an empty product equals $1$. Such a $\PP_u$ exists,
because $L^u$ is a nonnegative local martingale, satisfying $\sup_{t\leq T}
L^u_t\leq e^{A_T} C^{N_T}$ by \eqref{CC1}, and the latter variable is
integrable by \eqref{CC2}, so $L^u$ is indeed a uniformly integrable
martingale, with of course $\E(L^u_T)=1$. By Girsanov's theorem for point
processes, the predictable compensator of the measure $\mu$ under $\PP_u$ is
\[
\nu^u(dt,dx)=r(t,x,u_t) \nu(dt, dx)=r(t,x,u_t) \phi_t(dx) \,dA_t.
\]
We finally define the cost associated to every $u(\cdot)\in\aaa$ as
\[
J\bigl(u(\cdot)\bigr)=\E_u \biggl(\int_0^Tl(t,u_t) \,dA_t + \xi  \biggr),
\]
where $\E_u$ denotes the expectation under $\PP_u$.

Observe that, if $V_t=\int_0^t \sup_{u\in U}|l(s,u)| \,dA_s$, we have
\[
\E_u \biggl(\int_0^T\bigl|l(t,u_t)\bigr| \,dA_t \biggr)
\le\E_u \biggl( \int_0^T\sup_{u\in U}\bigl|l(t,u)\bigr| \,dA_t \biggr)
=\E\bigl(L^u_T V_T\bigr).
\]
Since $L^u$ is a nonnegative martingale and $V$ is continuous, adapted and
increasing, we deduce
\begin{eqnarray}
\label{CC4}&&\quad \E\bigl(L^u_T V_T\bigr)=\E
\biggl( \int_{0}^T L^u_t
\,dV_t \biggr)\le \E \biggl(\int_{0}^T
e^{A_t} C^{N_t}\sup_{u\in U}\bigl|l(t,u)\bigr|
\,dA_t \biggr)<\infty
\end{eqnarray}
by \eqref{CC3}.
Similarly,\vspace*{1pt} $\E_u(|\xi|)=\E(|\xi|L^u_T)\le \E (|\xi|e^{A_T} C^{N_T} )
<\infty$, and we conclude that under (\ref{CC3}) the cost $J(u(\cdot))$
is finite for every admissible control.

\begin{rema}  Suppose that the cost functional has the form
\[
J^1\bigl(u(\cdot)\bigr)=\E_u \biggl(\sum_{n\ge1 : S_n\le T}c(S_n,X_{n},u_{S_n}) \biggr)
\]
for some given predictable function $c$ on $\Omega\times[0,T]\times E\times U$
which is, for instance, nonnegative. By a standard procedure, we can reduce this control
problem to the previous one because
\begin{eqnarray*}
J^1\bigl(u(\cdot)\bigr)&=& \E_u \biggl(\int_0^T\!\int_E  c(t,x, u_{t}) \mu(dt,dx) \biggr)\\
&=&
\E_u \biggl(\int_0^T\!\int_E  c(t,x, u_{t})r(t,x, u_{t}) \phi_t(dx) \,dA_t \biggr).
\end{eqnarray*}
Thus, $J^1(u(\cdot))$ has the same form as $J(u(\cdot))$, with $\xi=0$
and with the function $l$ replaced by $l^1(t,u)=\int_E c(t,x, u)r(t,x,u)
 \phi_t(dx)$, so our forthcoming results can be applied.

Similar considerations obviously hold for cost functionals of the
form $J(u(\cdot))+J^1(u(\cdot))$.
\end{rema}

The control problem consists in minimizing $J(u(\cdot))$ over
$u(\cdot)\in\aaa$, and to this end a basic role is played by  the BSDE
\begin{eqnarray}
&& Y_t+\int_{(t,T]}\int
_EZ(s,x) \mu(ds,dx) = \xi +\int_{(t,T]} f
\bigl(s,Z_s(\cdot)\bigr) \,dA_s, \nonumber
\\[-8pt]
\label{CC5}
\\[-8pt]
\eqntext{t\in [0,T],}
\end{eqnarray}
with terminal condition $\xi$ being the terminal cost above, and with the
generator $f$ being the Hamiltonian function defined below. This is equation
\eqref{S11}, with $f$ only depending on $(\om,t,\ze)$, and indeed it comes
from an equation of type II via the transformation \eqref{S14}.

The Hamiltonian function $f$ is defined on $\Om\times[0,T]\times\ba(E)$ as
\begin{eqnarray}
&&\label{CC6}
f(\om,t,\ze)=\cases{ \displaystyle\inf_{u\in U}
\biggl(l(\om,t,u)+ \int_E\ze(x) r(\om,t,x,u)
\phi_t(\om,dx) \biggr), \vspace*{3pt} \cr
\qquad\mbox{if }\displaystyle\int_E\bigl|
\ze(x)\bigr| \phi_{\om,t}(dx)<\infty,\vspace*{3pt}
\cr
0,\qquad\mbox{otherwise}.}
\end{eqnarray}

We will assume that the infimum is in fact achieved, possibly at
many points. Moreover, we need to verify that the generator of the
BSDE satisfies the conditions required in the previous section, in particular
the measurability property, as expressed in \eqref{S12}, which does not
follow from its definition.
An appropriate assumption is the following one, since we
will see below in Proposition~\ref{PCC1} that it can
be verified under quite general conditions.

\renewcommand{\theass}{(C)}
\begin{ass}\label{cc}
For every predictable function $Z$ on $\Omega\times [0,T]\times E$
there exists a $U$-valued predictable process
(i.e., an admissible control) $\underline{u}^Z$ such that,
$dA_t(\omega)\PP(d\omega)$-almost surely,
\begin{eqnarray}
&& f\bigl(\om,t,Z_{\om,t}(\cdot)\bigr)
\nonumber
\\[-9pt]
\label{CC7}
\\[-9pt]
\nonumber
&&\qquad=l\bigl(\om,t,
\underline{u}^Z(\omega,t )\bigr)+ \int_E
Z_{\om,t}(x) r\bigl(\om,t,x,\underline{u}^Z(\om,t)\bigr)
\phi_t(\om,dx).
\end{eqnarray}
\end{ass}

Now, it is easy to check that all the required assumptions for the solvability
of the BSDE \eqref{CC5} are satisfied. Namely, using \eqref{CC1},
one easily proves the\vspace*{-2pt} inequality
\[
\bigl|f(\om,t,x,\ze)-f(\om,t,x,\ze')\bigr|\le C\int_E \bigl|\ze(y) -\ze'(y)\bigr|
\phi_{\om,t}(dy),
\]
whereas $f(\om,t,0)=\inf_{u\in U} l(\om,t,u)$. Then, in view of \eqref{CC3},
we see that \eqref{S12} and \eqref{S23} are satisfied, with $L=C$ and
$L'=0$, hence $\be>1+L+L'$ and $\al>L$. We thus conclude from Theorem~\ref{TS2} that the BSDE has a unique solution
$(Y,Z)\in \laa^1_{\alpha,\beta}$.
The corresponding admissible control $\underline{u}^Z$, whose
existence is required in Assumption \ref{bb}, will be denoted as $u^*$.

\begin{theo}\label{TCC1}
Assume \ref{aa}, \ref{bb} and \ref{cc}. Then, with $(Y,Z)$ and $u^*$
as above, the admissible control $u^*(\cdot)$ is optimal, and
$Y_0= J(u^*(\cdot))=\break \inf_{u(\cdot)\in\aaa }J(u(\cdot))$ is the
minimal cost.
\end{theo}

\begin{pf}
Fix $u(\cdot)\in\aaa$. We first show that
$\E_u \int_{0}^T\!\int_E |Z(t,x)| \nu^u(dt,dx)<\infty$. Indeed, setting $V_t=
\int_0^t\int_E  |Z(s,x)| r( s,x,u_s)  \nu(ds,dx)$ and arguing as in
(\ref{CC4}),
\begin{eqnarray*}
&& \E_u \biggl(\int_{0}^T \!\int_E |Z(t,x)| \nu^u(dt,dx) \biggr)\\[-2pt]
&&\qquad=
\E_u \biggl(\int_{0}^T \!\int_E  \bigl|Z(t,x)\bigr|  r( t,x,u_t)  \nu(dt,dx) \biggr)\\[-2pt]
&&\qquad=\E \bigl(L^u_T V_T\bigr)=\E \biggl(\int_{0}^T L^u_t \, dV_t \biggr)\le
\E \biggl(\int_{0}^T e^{A_t} C^{N_t}  \,dV_t \biggr)\\[-2pt]
&&\qquad= \E \biggl(\int_{0}^T\! \int_E  e^{A_t} C^{N_t} \bigl|Z(t,x)\bigr| r( t,x,u_t)
\nu(dt,dx) \biggr)\\[-2pt]
&&\qquad \le
 C \E \biggl(\int_{0}^T \!\int_Ee^{\be A_t} \al^{N_t} \bigl|Z_t(x)\bigr| \nu(dt,dx) \biggr),
\end{eqnarray*}
which is finite, since $(Y,Z)\in \laa^1_{\alpha,\beta}$.
By similar arguments, we also check that
\begin{eqnarray*}
&&\E_u \biggl(\int_{0}^T \bigl|f\bigl(t,Z_t(\cdot)\bigr)\bigr| \,dA_t \biggr)
\\[-2pt]
&&\qquad=\E \biggl(\int_{0}^T L_t^u  \bigl|f\bigl(t,Z_t(\cdot)\bigr)\bigr| \,dA_t \biggr)\le
\E \biggl(\int_{0}^Te^{A_t} C^{N_t} \bigl|f(t,Z_t\bigl(\cdot)\bigr)\bigr| \,dA_t \biggr) \\[-2pt]
&&\qquad\le \E \biggl(\int_{0}^T e^{ A_t} C^{N_t}
 \biggl(C\int_E \bigl|Z(t,x)\bigr| \phi_t(dx)+\bigl|f(t,0)\bigr| \biggr) \,dA_t \biggr)<\infty.
\end{eqnarray*}
Setting $t=0$ and taking the $\PP_u$-expectation in the BSDE (\ref{CC5})
we therefore obtain
\[
Y_0+\E_u \biggl(\int_{0}^T \!\int_EZ(t,x)  r(t,x,u_t)  \nu(dt,dx) \biggr)=
\E_u(\xi)+\E_u \biggl(\int_0^T f\bigl(t,Z_t(\cdot)\bigr) \,dA_t \biggr).
\]
Adding $\E_u (\int_{0}^T l(t,u_t) \,dA_t )$
to both sides, we finally obtain the equality
\begin{eqnarray*}
&& Y_0+\E_u \biggl(\int_{0}^T  \biggl(l(t,u_t)+\int_EZ(t,x) r(t,x,u_t) \phi_t(dx)
 \biggr) \,dA_t \biggr)\\
&&\qquad= J\bigl(u(\cdot)\bigr)+\E_u \biggl(\int_{0}^T f\bigl(t,Z_t(\cdot)\bigr) \,dA_t \biggr)\\
&&\qquad=J\bigl(u(\cdot)\bigr)+\E_u \biggl(\int_{0}^T\inf_{u\in U} \biggl(l(t,u)+\int_E
Z(t,x) r(t,x,u_t),\phi_t(dx) \biggr) \,dA_t \biggr).
\end{eqnarray*}
This implies
immediately the inequality $Y_0 \le  J(u(\cdot))$
for every admissible control, with an equality if $u(\cdot)=u^*(\cdot)$.
\end{pf}

Assumption \ref{cc} can be verified in specific
situations when it is possible to compute explicitly the function
$\underline{u}^Z$. General conditions for its validity can also be
formulated using appropriate measurable selection theorems, as in  the
following proposition.

\begin{prop} \label{PCC1} Suppose that $U$ is a compact metric space with its
Borel $\sigma$-field $\ua$ and that the functions
$r(\omega,t,x,\cdot),l(\omega,t,\cdot)$  are continuous on $U$ for every
$(\om,t,x)$. Then if further \eqref{CC1} holds, Assumption \ref{cc} is satisfied.
\end{prop}

\begin{pf}
For every predictable function $Z$ set
$G^Z= \{(\om,t) : \int_E|Z(\omega,t,\break x)| \phi_{\om,t}(dx)=\infty \}$
and define a map $F^Z:\Omega\times [0,T]\times U\to\R$ by
\begin{eqnarray*}
&& F^Z(\omega,t,u)\\
&&\qquad=\cases{
\displaystyle l(\omega,t,u)+\int_E Z(\om,t,x)   r(\om,t,x,u)  \phi_t(\om,dx),  &\quad$\mbox{if }
(\om,t)\notin G^Z$,\vspace*{3pt}\cr
0,   &\quad $\mbox{if }(\om,t)\in G^Z$.}
\end{eqnarray*}
Then $F^Z(\omega, t,\cdot)$ is continuous for every $(\om,t)$
and $F^Z$ is a predictable function on $\Om\times[0,T]\times U$. By a
classical selection theorem (see, e.g., Theorems 8.1.3 and 8.2.11 in \cite{AF}
there exists a $U$-valued function $\underline{u}^Z$ on $\Omega\times [0,T]$
such that $F^Z(\omega,t,\underline{u}^Z(\omega,t)) =
\inf_{u\in U} F^Z(\omega,t,u)$ for every $(\omega,t)\in \Omega\times [0,T]$
[so that~\eqref{CC7} holds true for every $(\omega,t)$] and such that
$\underline{u}^Z$ is measurable with respect to the completion of the
predictable $\sigma$-algebra in $\Omega\times [0,T]$ with respect to the
measure $dA_t(\omega)\PP(d\omega)$. After modification on a null set, the
function $u^Z$ can be made
predictable, and \eqref{CC7} still holds, as it is understood as an equality
for $dA_t(\omega)\PP(d\omega)$-almost all $(\omega,t)$.
\end{pf}


\printaddresses
\end{document}